\documentclass[12pt,a4paper,english]{smfart} 

\usepackage{fullpage}

\usepackage[table]{xcolor}

\usepackage[utf8]{inputenc}

\usepackage[OT2,T1]{fontenc}

\usepackage{amsfonts,lmodern}
\usepackage{amssymb,mathabx}
\usepackage{ifthen}
\usepackage{pstricks}
\usepackage{latexsym,eurosym,ulem}
\usepackage[french]{babel} 
\usepackage{footmisc}

\usepackage{alltt}
\usepackage{calrsfs}
\usepackage{fancyhdr}
\usepackage{graphicx}

\usepackage{amssymb}
\usepackage{amsmath}
\usepackage{smfthm,mathabx}
\usepackage{enumerate}
\usepackage{amsthm}
\usepackage{amsfonts}
\usepackage{graphicx}
\usepackage[colorlinks=true,linkcolor=black,citecolor=black,urlcolor=black]{hyperref}
\usepackage{fancyhdr}
\usepackage[all,cmtip]{xy}
\usepackage{alltt}
\usepackage{footmisc}
\usepackage{smfthm}

\font\rurm=wncyr10 scaled \magstep1

\def\sha{{{\textnormal{\rurm{Sh}}}}}

\def\Sha{{\sha}^2}

\def\Z{{\set Z}}

\def\R{{\mathbb{R}}}

\def\Z{{   \mathbb Z }}
\def\Q{{\mathbb Q}}

\def\F{{\mathbb F}}

\def\L{{\mathcal L}}

\def\E{{\mathcal E}}

\def\m{{\mathfrak m}}

\def\q{{\mathfrak q}}
\def\wp{{\q}}

\def\gl{{\mathfrak{gl}}}
\def\sl{{\mathfrak{sl}}}

\def\gr{\mathfrak{gr}}

\def\Matrix#1#2#3#4{{\left(\begin{array}{cc} #1 & #2 \\ #3& #4 \end{array}\right)}}

\def\sha{{{\textnormal{\rurm{Sh}}}}}
\def\CyB{{{\textnormal{\rurm{B}}}}}
\def\Ktame{K^{{\rm ta}}}
\def\Gtame{G^{{\rm ta}}_K}
\def\GtameS2{G_{S_2}}
\def\Sha{{\sha}^2}
\def\Gal{\mathrm{Gal}}
\def\Res{\mathrm{Res}}
\def\Inf{\mathrm{Inf}}
\def\SL{\mathrm{SL}}
\def\GL{\mathrm{GL}}
\def\Rad{\mathrm{Rad}}

\def\Adx{\mathrm{ad}_x}
\def\Ady{\mathrm{ad}_y}
\def\Ad{\mathrm{ad}}

\def\brackxy{[x\,y]}
\def\brackxz{[x\,z]}
\def\brackyz{[y\,z]}

\def\brackLL{[L\,L]}
\def\im{\mathrm{Im}}

\newtheorem{Theorem}{Theorem}

\newtheorem*{Question}{Question}

\newtheorem{Corollary}[Theorem]{Corollary}

\begin{document}


\title{On tamely ramified infinite Galois extensions}
\author{Farshid Hajir, Michael Larsen, Christian Maire, Ravi Ramakrishna}
\address{Department of Mathematics \& Statistics, University of Massachusetts, Amherst, MA 01003, USA}
\address{ Department of Mathematics, Indiana University, Bloomington, IN 47405, USA}
 \address{FEMTO-ST Institute, Universit\'e de Franche-Comt\'e, CNRS,  15B avenue des Montboucons, 25000 Besan\c con, FRANCE} 
\address{Department of Mathematics, Cornell University, Ithaca, NY 14853-4201 USA}
\email{hajir@math.umass.edu, mjlarsen@indiana.edu, christian.maire@univ-fcomte.fr, ravi@math.cornell.edu}

\begin{abstract}
For a number field $K$, 
we consider $\Ktame$ the maximal tamely ramified algebraic extension of~$K$, and its Galois group $\Gtame= \Gal(\Ktame/K)$. Choose a prime $p$ such that $\mu_p \not \subset K$. Our guiding aim is to characterize the finitely generated pro-$p$ quotients of~$\Gtame$.
We give a {unified point of view} by introducing  the notion of  {\it stably inertially generated}  pro-$p$ groups~$G$, for which linear groups  are archetypes.  
This key notion {is compatible} with local {\it tame liftings} as used in the Scholz-Reichardt Theorem. We realize every finitely generated pro-$p$ group~$G$ which is stably  inertially generated as a quotient of $\Gtame$. 
Further examples of groups that we realize as quotients of $\Gtame$ include 
congruence subgroups of special linear groups over $\Z_p\ldbrack T_1,\cdots, T_n \rdbrack$. Finally, we give classes of groups which cannot be realized as quotients of  $G^{\rm ta}_{\Q}$. 
\end{abstract}

\thanks{The second author was partially supported by the National Science Foundation  grant DMS-2001349.
The third author  was partially supported by  the EIPHI Graduate School (ANR-17-EURE-0002). The fourth author was partially supported by Simons Collaboration grant \#524863. He also thanks FEMTO-ST for its hospitality and wonderful research environment during his visit there in the spring of 2022. The first, third and fourth authors were supported by ICERM for a Research in Pairs visit in January, 2022 and by the Western Academy for Advanced Research (WAFAR) for a visit during the summer of 2023.}

\date{\today}




\maketitle

Let $K$ a number field.
Consider $\Ktame$, the maximal tamely ramified algebraic extension of~$K$.
Let $\Gtame= \Gal(\Ktame/K)$.  Choose a prime $p$. 
The maximal  pro-$p$ quotient of $\Gtame$ has an infinite number of generators. We are interested in the following question:

\begin{Question} What are the finitely generated pro-$p$ quotients of $\Gtame$?
\end{Question}

\smallskip

We make some first observations  regarding these possible quotients. 

$-$    
By a classical theorem of Scholz-Reichardt, for $p$ odd, every finite $p$-group $G$ is a quotient of 
$G^{\rm ta}_\Q$. 
See Serre \cite[Chapter 2, \S 2.1]{Serre} for a treatment which in fact generalizes their result. 
Shafarevich resolved the case of $p=2$ by introducing the shrinking process; see also a recent work of Schmid \cite{Schmid}.\color{black}

$-$  Recall that a pro-$p$ group is called FAb if its open subgroups have finite abelianization. Finitely generated quotients of $\Gtame$ are FAb by class field theory.

$-$  
The previous  constraint is related to   
the tame  Fontaine-Mazur conjecture \cite[Conjecture 5a]{FM} which implies that every finitely ramified   $p$-adic analytic quotient of $\Gtame$  must be finite. By finitely ramified we mean the set of ramification is finite, though
the ramification index at a particular prime could be infinite. The only established case of the tame Fontaine-Mazur conjecture is that of 1-dimensional representations (which amounts to the FAb property of finite generated quotients of $\Gtame$), and this is a good measure of both the importance and the difficulty of studying the finitely generated quotients of $\Gtame$. Much of this difficulty lies in the subtleties of the tame Galois cohomology - in the wild analogue (where all primes above $p$ are allowed to ramify), Poitou-Tate duality renders many of the calculations much more straightforward.

$-$ 
In \cite{RR2} {surjective} representations $G_{\mathbb Q}  \twoheadrightarrow \SL_2({\mathbb Z}_p)$ 
are constructed that are unramified  at $p$ provided the mod $p$ reduction is unramified  at $p$.  
These representations are almost certainly ramified at infinitely many tame primes as the Fontaine-Mazur conjecture  predicts, though proving so definitively seems very difficult for even representations. 

\medskip

Choose $p$ such that    $\mu_p \not \subset K$, where $\mu_p$ denotes the  $p$th roots of unity. We make this hypothesis to guarantee linear disjointness of various field extensions which in turn allows us to apply Chebotarev's theorem and successively resolve obstruction problems. See Proposition \ref{proposition_lifting}. When 
$\mu_p \subset K$, these issues are technically much more involved and it is not clear that for a fixed number field $K$, the statements we prove for large enough $p$ always hold for $p$ such that $\mu_p \subset K$.
As $\mu_2$ is in every number field, we require
$p  >2$.
Under the assumption that $K$ does not contain a primitive $p$th root of unity,  we use embedding techniques to characterize many finitely generated pro-$p$ groups~$G$ that can be realized as quotients  of $\Gtame$. These constructions are explained below in Theorem~\ref{TheoremA} and its three corollaries C, D, and E. We also give classes of groups that are not quotients of $\Gtame$; see Theorem~\ref{NoToralQuotient}. 
To obtain these results, we give a {unified point of view} by introducing  the notion of  
{\it stably inertially generated}  pro-$p$ groups $G$,
for which linear groups  are archetypes.  
This key notion has a certain compatibility with local {\it tame  liftings} as used in the Scholz-Reichardt Theorem  \cite[Chapter 2, \S 2.1]{Serre} extended by Neukirch \cite{Neukirch}.

\smallskip

The tame $p$-adic Lie Galois extensions which we construct have the further property that  the set of primes that split completely is infinite. This phenomenon arises  in the wild ``{finitely ramified case}'' when the $p$-adic Lie group is nilpotent as observed in 
\cite{HM_cambridge}, Proposition 4.1. In that situation the set of primes that split completely is related to the $\mu$-invariant of the relevant $p$-adic Lie extension.
Here we show that if  one allows infinitely many ramified primes, then infinitely many primes can split completely for the groups   $\SL_m^k(\Z_p)$ below. This result is similar to that of 
\cite{KLR} where the characteristic polynomials of almost all Frobenius elements are arranged to have pure algebraic roots.  It seems very difficult to achieve such results using only Galois cohomology without the infinite ramification hypothesis.

\medskip

All  generation statements in this paper refer to topological generators. 
We will always denote by $G$ a  finitely generated pro-$p$ group.
We call an element $y \neq 1$ of $G$  {\it inertial} if there exists $x \in G$ such that $$[x,y] = y^\lambda$$  for some $0 \neq \lambda \in \Z_p$.
We say a pro-$p$ group $G$ is \textit{inertially generated} if it is  generated by 
inertial elements.
Indeed, this definition is motivated by the standard relation
$[\sigma,\tau]=\tau^{q-1}$ for the Galois group over $\Q_q$ of its maximal pro-$p$ extension when 
$q \cong 1$ mod $p$. 

We say the pro-$p$ group $G$ is  
{\it stably inertially generated}  if each group $P_n(G)$ of  the $p$-central series $(P_n(G))_n$ of~$G$ is inertially generated. 
See \S \ref{section_locally_inertially} for the definitions and examples.

Our main theorem is:

\begin{Theorem}\label{TheoremA}
      Let $G$ be a finitely generated  pro-$p$ group that is 
 stably   inertially generated. Then there exists a Galois extension $L/K$ in $K^{\rm ta}/K$ such that  $\Gal(L/K)\cong G$. Moreover $L$ can be taken such that the set of primes of $K$ that split completely in $L$
 is infinite.
     \end{Theorem}

      Let us say few words regarding our strategy, which   is classical.  We filter our pro-$p$ group $G$ with subgroups $H_n$ such that
\begin{itemize}
\item $G/H_2\cong (\Z/p)^d$ where $d$ is the minimal number of generators of $G$ and
\item  $H_n/H_{n+1} \cong \Z/p$ for $n \geq 2$.
\end{itemize}
We then construct $\Gtame \twoheadrightarrow G/H_2$ as the base case of our induction and inductively
build compatible homomorphisms $\Gtame \twoheadrightarrow G/H_n$. Taking the limit solves the problem. 
It is easy to arrange that all obstructions are local, but then one needs to remove the local obstructions. At each stage of the induction we must allow ramification at another prime {\it and} make certain there is no local 
obstruction to lifting to the next step at this new prime. This requires a Galois cohomological argument and that
$G$ is stably inertially generated.

     \medskip
 
We say $G$ is \textit{torsion-generated} if it is generated by elements of finite order.
We say it is \textit{stably torsion-generated} if the subgroups $P_n(G)$ are torsion-generated for all $n\geq 1$.

     \begin{Corollary} \label{Unstable torsion}
      Let $G$ be a stably  torsion-generated pro-$p$ group.
     Then  $G$ can be realized as a quotient $\Gal(L/K)$ of $\Gtame$ and the extension $L$ can be taken such that 
    the set of primes of $K$ that split completely in $L$
 is infinite.
            \end{Corollary}

This corollary can be compared to Serre's pro-$p$ version of the Scholz-Reichardt Theorem \cite[Theorem 2.1.11]{Serre}.  The latter theorem does not require $G$ to be finitely generated, but it assumes that $G$ has finite exponent, so, in particular, its elements are all of finite order. By a famous result of Zelmanov, any finitely generated pro-$p$ group $G$ of finite exponent is itself finite, so the pro-$p$ version of Scholz-Reichardt becomes the classical result in the finitely generated case.
 
Note that stably torsion-generated groups 
which have non-torsion elements exist. 
 Consider the free group $F_2$ on two generators $x$ and $y$. Let $\{a_i\}^\infty_{i=1}$ 
 be an enumeration of the countably many elements of this group. Impose the relations $R=\{a^{p^{m_i}}_i\}$ with $m_i \to \infty$ in such a way that the Golod-Shafarevich power series is negative on $(0,1)$ (see Theorem 7.20 of \cite{Koch}). This implies the pro-$p$ completion $H$ of $F_2/\langle R\rangle$  is infinite. Since the words of $F_2$ are dense in $H$, 
By a theorem of Zelmanov \cite{Ze}, $H$ contains a free pro-$p$ group as a subgroup and is therefore not torsion.

\medskip

Assume $p>2$. For a complete Noetherian local ring $A$ with maximal ideal $\m$ such that $A/\m\cong\F_p$, 
and $k\geq 1$, set 
$$\SL_m^k(A):=\ker\left(\SL_m(A)\rightarrow \SL_m(A/\m^k)\right).$$
Then,
$$\SL_m^k(A)/\SL_m^{k+1}(A) \cong (\m^k/\m^{k+1})\otimes_{\Z} M_m^0(\Z),$$
where $M_m^0(\Z)$ denotes $m\times m$ integer matrices with trace $0$.  As $A$ is Noetherian this is a finite dimensional vector space over $\F_p$, so
$$\widehat{\SL}_m^1(A) := \projlim_k \SL_m^1(A)/\SL_m^k(A)$$ 
is a pro-$p$ group  generated by any set of elements in $\SL_m^1(A)$ 
lifting an $\F_p$-basis of
 $\SL_m^1(A)/\SL_m^2(A) \cong (\m/\m^2)\otimes_\Z M_m^0(\Z)$. 
The natural homomorphism $\SL_m^1(A)\to \widehat{SL}_m^1(A)$ is injective by the Krull intersection theorem and surjective by the completeness of $A$.
Therefore, $\SL_m^k(A)$ is a   finitely generated    pro-$p$ group for all $k$.

\begin{Corollary} \label{corollary_padic}  For $k,m,n \geq 1$,  
 the groups $\SL_m^k(\Z_p\ldbrack T_1,\cdots, T_n \rdbrack)$ are quotients of~$\Gtame$.
 In particular each $\SL_m^k(\Z_p)$ is a quotient of $\Gtame$. Moreover these quotients can be chosen to correspond to a Galois extension $L/K$ in which infinitely many primes split completely.
\end{Corollary}

This result {is (in spirit) an extension of} those of \cite{RR1}, and \cite{KLR} for $\SL_2(\Z_p)$ and \cite{RR2} for $\SL_2(\Z_p\ldbrack T_1,\cdots, T_n \rdbrack)$. 

\medskip

The notion of 
stably inertially generated pro-$p$ group is particularly well-adapted for $p$-adic Lie groups.
To each $p$-adic analytic group $G$ one can attach in a natural way a ${\mathbb Q}_p$-Lie algebra $L_G$ \cite{Lazard}, \cite{Dixon};  we  recall  this principle in \S \ref{section_dictionary}. 
We   focus on  the case that $L_G$ is simple. There are only two possibilities: $L_G$ is pluperfect or toral.
We say a Lie algebra $L$ is {\it toral} if $\Adx$ is semisimple for all $x\in L$, and it is {\it pluperfect} if it admits no non-trivial toral quotient algebra.
See \S \ref{section_toral_pluperfect} for the definitions  and further discussion of these concepts. Using Theorem \ref{TheoremA}
 we prove:

\begin{Corollary}\label{SimplePluperfect} Let $G$ be a $p$-adic analytic group with
pluperfect Lie algebra~$L_G$. Then there exist 
 homomorphisms \color{black} $\rho:\Gtame \to G$ such that the image  of $\rho$ is open in $G$. 
\end{Corollary}

In the other direction, by class field theory one knows that uniform  abelian  quotients of~$\Gtame$ are trivial. See \S \ref{section_padic} for the definition of a uniform group. In fact  uniform abelian  groups are a special case of a more general family,  uniform toral groups (the pro-$p$ group is uniform and its Lie algebra is toral). 
Examples of uniform toral FAb  groups are  those whose Lie algebra  is the set of  trace zero elements of a skew field. See Section \ref{subsection_local}.

\begin{Theorem}\label{NoToralQuotient} If the $p$-class field tower of $K$ is finite then
the pro-$p$ group $\Gtame$ has no non-trivial uniform toral  quotient.  In particular, for a fixed $K$, the latter statement is true for all large enough primes $p$.
\end{Theorem}
Uniform toral  groups are finitely generated. Given a number field $K$, one expects by Gras' conjecture \cite{Gras-CJM} that the Galois group over $K$ of its maximal pro-$p$ extension unramified outside primes above $p$ is free pro-$p$ for $p$ large enough. The number of generators would be $r_2(K)+1$ where $r_2(K)$ is the number of pairs of complex embeddings of $K$.
The uniform toral FAb groups of Theorem~\ref{NoToralQuotient} are quotients of such groups for large enough $r_2(K)$, so they provide (conjectural)  examples of wildly ramified Galois groups that cannot arise as quotients of $\Gtame$.
The Fontaine-Mazur conjecture
for tame extensions predicts 
Theorem~\ref{NoToralQuotient}  holds without the class field tower hypothesis.

In order to 
refine Corollary~\ref{SimplePluperfect} 
consider the following context.

Let $T$ be a (possibly infinite) set of primes of $K$, 
 $K^{{\rm ta},T}/K$ be the maximal $T$-split  extension of $K$ in $\Ktame$, and
$G^{{\rm ta},T}_K=\Gal(K^{{\rm ta},T}/K)$. 

Assuming the condition that $\mu_p \not \subset K$, all the previous results can be extended to $G^{{\rm ta},T}_K$ without difficulty.

\medskip

Set 
  $$\displaystyle{\alpha_T:= \sum_{\wp\in T} \frac{\log N(\wp)}{N(\wp)-1}} \mbox{ and } \displaystyle{\alpha^{\rm GRH}_T:= \sum_{\wp\in T} \frac{\log N(\wp)}{\sqrt{N(\wp)}-1}},$$
  where $N(\wp):=\# {\mathcal O}_K/{\mathfrak q}$.
 
For  $v|\infty$ set
 $$ \alpha_\q= \left\{ \begin{array}{cl} \frac{1}{2}(\gamma+ \log4\pi) &  v \mbox{ is real} \\
\gamma+\log 2\pi & v  \mbox{ is complex} \end{array}\right. \mbox{ and }
 \alpha^{\rm GRH}_\q= \left\{ \begin{array}{cl} \frac{1}{2}(\frac{\pi}{2}+\gamma + \log8\pi ) &  v \mbox{ is real} \\
\gamma+ \log 8\pi& v  \mbox{ is complex} \end{array}\right.$$ 
{where $\gamma$ is Euler's constant.}
Theorem~\ref{SplittingTheorem} below, unlike Theorem~\ref{NoToralQuotient}, requires no condition on the $p$-class field tower of $K$. 
\begin{Theorem}\label{SplittingTheorem} Let $d_K$ to be the absolute discriminant of $K$.\\
\noindent
1) Take $T$ such that $\alpha_T  + \sum_{v|\infty} \alpha_\q > \log \sqrt{|d_K|} $.  Then $G^{{\rm ta},T}_K$ has no non-trivial  uniform toral  quotient.\\
2) Assume the GRH and  $T$ such that $\alpha^{\rm GRH}_T  + \sum_{v|\infty} \alpha^{\rm GRH}_\q > \log \sqrt{|d_K|} $.  Then $G^{{\rm ta},T}_K$ has no non-trivial  uniform toral  quotient.
\end{Theorem}

 In \S~\ref{section:one} we establish the Galois cohomological machinery we need, Proposition~\ref{proposition_lifting} being the important  technical result.
In \S~\ref{section:two} we perform the inductive lifting process, making key use of Definition~\ref{xwhereweneedit} to
prove Theorem~\ref{TheoremA}.
In \S~\ref{section:three} we establish results on $p$-adic Lie algebras and classes of pro-$p$ groups that allow us to obtain the rest of our results.

\vskip1em
{\bf Notations.} 
Throughout this article $p>2$ is an odd prime number and $G$ is a {finitely generated} pro-$p$ group.

$\bullet$  Set $G^{ab}:=G/\overline{[G,G]}$, $G^{p,el}:=G^{ab}/(G^{ab})^p$, and $d(G):= \dim G^{p,el}$.

$\bullet$ Let $(P_n(G))$ be the $p$-central series of~$G$: $P_1(G)=G$ and for $n\geq 1$,  
$P_{n+1}(G)=\overline{P_n(G)^p[G,P_n(G)]}$.
The sequence $(P_n(G))$ forms a basis of open normal subgroups of $G$. 

$\bullet$ 
All cohomology groups have coefficients in $\Z/p$ with trivial action so
we write $H^i(G)$ for $H^i(G,\Z/p)$.

$\bullet$ $K$ is a number field with $\mu_p \not \subset K$ and $\Ktame$ is the maximal tamely ramified extension of $K$ and $\Gtame=\Gal(\Ktame/K)$.

$\bullet$ For a prime $\q$ of $K$ we denote the cardinality of ${\mathcal O}_K/{\q}$ by $N(\q)$.

$\bullet$ For a finite set of primes $S$ of residue characteristic different from $p$, set $K_S$ to be the maximal pro-$p$ extension of $K$ unramified outside $S$ and $G_S=\Gal(K_S/K)$.

$\bullet$ $K_\q$ is the completion of $K$ at the prime $\q$, $G_\q:=\Gal(\overline{K}_\q/K_\q)$ and $U_\q$ is the group of units of $K_\q$.  We let $\sigma_\q$ and $\tau_q$ be, respectively, a 
Frobenius element and a generator of inertia in the Galois group of the maximal pro-$p$ extension of $K_\q$.

$\bullet$ $J$ is the idele group of $K$. The subgroup $U=\prod_\q U_\q\subset J$ are those ideles that are locally units everywhere and $U_S \subset U$
are those ideles in $U$ that are $1$ at every $\q\in S$.

$\bullet$ $V_S:=\{x\in K^\times | x\in K^{\times p}_\q , \ \forall \q \in S; x \in U_\q K^{\times p}_\q\ \forall \q \}$ and
$\CyB_{S} = (V_{S}/K^{\times p})^\wedge$   is called the Selmer group.  

$\bullet$ For $S \supseteq T$, set $V_T^S=\{x\in K^\times \mid  \ x\in K ^{\times p}_\q   \,\, \forall \,\, \q \in T 
\mbox{ and } x\in U_\q  K^{\times p}_\q  \,\,  \forall  \,\, \q \notin S   \}$
and  $\CyB^S_T = (V^S_T/ K^{\times p})^\wedge.$

\section{The local-global principle}\label{section:one}

\subsection{The Shafarevich and Selmer groups}

Let $Z$ be a finite set of tame primes and recall $G_{Z}$ is the Galois group of the maximal extension of $K$ unramified outside primes above $Z$. 
 We will analyze our obstructions in group cohomology via the exact restriction sequence
\begin{eqnarray}\label{local_global} 0\to \Sha_{Z}\to H^2(G_{Z}) \stackrel{\Res_{Z}}\longrightarrow \prod_{\q \in Z} H^2(G_\q) \end{eqnarray}
where   $\Sha_{Z}$ is defined as the kernel of  the restriction map  $\Res_Z$. 
The groups $V_Z$ and $\CyB_Z$ from the Notations are crucial to the study of $G_{Z}$.
We record two important  results   (parts $(i)$ and $(iv)$) and several basic facts about these groups. 

\begin{lemm}\label{lemm:Shaprops} Let $K$ be a number field. Then \\
(i) $\Sha_{S} \hookrightarrow \CyB_{S}$.\\
(ii) $S_1 \subset S_2 \implies V_{S_2} \subseteq V_{S_1}$.\\
(iii) Let $U_K$ be the units  of the number field  $K$ and $Cl_K[p]$ be the $p$-torsion of the class group of $K$. There is an exact sequence
$$1 \to U_K/U^p_K \to V_\emptyset/K^{\times p} \to Cl_K[p] \to 1.$$
As $\mu_p \not \subset K$ we have $\dim \CyB_{K,\emptyset} =  r_1(K)+r_2(K) -1 +\dim Cl_K[p]$.  \\
(iv)  There is a finite set $Z_0$ of 
$\dim \CyB_{K,\emptyset}$ \color{black} tame primes such that
$\CyB_{Z_0}= 0$  so $ \Sha_{S}=0$ for any $S\supseteq Z_0$. \\
(v) For $S\supseteq T$, 
 $V^{S\cup T}_T \subseteq V^S_\emptyset$.\\
(vi) $V_T = U_TJ^p \cap K^\times$ so $V_T/K^{\times p} = (U_TJ^p \cap K^\times)/K^{\times p}\cong (U_TJ^p \cap J^p K ^\times)/J^p$.\\
(vii) $V_T^S = (U_TJ^p  \prod_{\q \in S\setminus T} K_\q^\times) \cap K^\times$ so 
 $V^S_T/K^{\times p} = (U_TJ^p  \prod_{\q \in S\setminus T} K_\q^\times \cap J^pK^\times)/J^p $.
\end{lemm}
\begin{proof}
$(i)$ See  \cite[Chapter 11, Theorem 11.3]{Koch}. \\
$(ii)$ This is immediate from the definition.\\
$(iii)$  One can alternatively define $V_\emptyset = \{ \alpha \in K^\times | (\alpha)=I^p\}$. 
Let $I$ be an ideal representing a class of $C$ of $CL_K[p]$ so $I^p=(\alpha)$. The map $\alpha \mapsto C$ is surjective and one easily sees the kernel is $U_K/U^p_K$. 
\\
$(iv)$
This is consequence of  the finiteness of $\CyB_{K,\emptyset}$: see for example  Theorem 1.12 of \cite{HMRI}.\\ 
$(v)$, $(vi)$ and $(vii)$ These are immediate from the definitions.
\end{proof}

\begin{rema}
It is worth remarking the injection $(i)$ above is an isomorphism if $S$ contains all primes above $p$. The failure
of this map to be surjective in the tame setting can be thought of as the reason tame Galois cohomology is difficult.
\end{rema}

\subsection{Global cohomology classes with given local conditions}

\subsubsection{} 

Proposition~\ref{proposition_10.7.9} below is Proposition 10.7.9 of \cite{NSW}. 
The proof in \cite{NSW} invokes previous results from that text using Poitou-Tate duality.
We give a proof more in the spirit of
\S 11.3 of \cite{Koch}, using only class field theory.
 \begin{prop} \label{proposition_10.7.9} Let $S \supseteq T$ be finite sets of primes. 
 There is an exact sequence
 $$0 \to H^1(G^T_S,\Z/p) \to H^1(G_S,\Z/p) \to \prod_{\q \in T} H^1(G_\q,\Z/p) \to  \CyB^S_{S\setminus T} 
 \to \CyB_S \to 0$$
 where $G_S$ is the Galois group of the maximal extension of $K$ unramified outside $S$ and 
 $G^T_S$ is the Galois group of the maximal extension of $K$ unramified outside $S$ split completely at $T$.
    \end{prop}
    \begin{proof} 
    Recall $U_X\subset J$ are those ideles that are $1$ at $\q \in X$ and units at $\q \notin X$. \color{black}
    Consider the sequence
$$
0 \to \frac{U_SJ^p \cap J^pK^\times}{J^p}\stackrel{\psi_1}{\to} 
\frac{(U_{S\setminus T} J^p  \prod_{\q \in T} K^\times_\q) \cap J^pK^\times}{J^p}  
\stackrel{\psi_2}{\to} \prod_{\q \in T} \frac{K^\times_\q}{K^{\times p}_\q} $$
$$\stackrel{\psi_3}{\to}
\frac{J}{U_SJ^pK^\times} \stackrel{\psi_4}{\to} \frac{J}{(U_SJ^p  \prod_{\q \in T} K^\times_\q)K^\times} \to 0
$$
where $\psi_1$ is the natural inclusion, $\psi_2$ is the projection to the $T$-coordinates, $\psi_3$ is extension from $T$-coordinates to the ideles by including the component $1$ at all $\q \notin T$ and taking the quotient by $U_SJ^pK^\times$, and $\psi_4$ is the natural projection.
We will prove this sequence is exact and dualizing will give the result.

\smallskip
Exactness at the first, fourth and fifth terms is clear. 
We now show the sequence is a complex. 
That $\psi_2\circ\psi_1$ is trivial follows from the fact that, since $S\supseteq T$,  it is trivial on $U_SJ^p$.
We show $\psi_3\circ \psi_2$ is trivial. Let $u_{S\setminus T} j_1^p t \in 
(U_{S\setminus T} J^p  \prod_{\q \in T} K^\times_\q) \cap J^pK^\times$ so we can write
$u_{S\setminus T} j_1^p t =j_2^p\gamma$ where $\gamma \in K^\times$ and thus
$u_{S\setminus T} j^p t =\gamma \in K^\times$ 
where $j=j_1/j_2$.
Up to $p$th powers of ideles, we see $\gamma$ is a unit outside of $S$, a $p$th power at $S \setminus T$  and 
$$\psi_3(\psi_2( u_{S \setminus T} j_1^p t)) = \psi_3(\psi_2( u_{S \setminus T} j^p t))
=\psi_3(\psi_2( \gamma))$$ which has component $\gamma$ at $\q \in T$ and $1$ elsewhere. 
Let $\tilde{u}_S$ be the idele with component $1$ at $\q \in S$ and $\gamma$ elsewhere. 
As $\gamma$ is, up to $p$th powers,  a unit outside of $S$ we see $\tilde{u}_S \in U_S J^p$ and $\gamma^{-1}\tilde{u}_S \psi_3(\psi_2( \gamma))$ has component $\gamma^{-1}$ at $\q \in S \setminus T$ and $1$ elsewhere. 
As $\gamma^{-1}$ is a $p$th power at $S \setminus T$, we have, in $J/(U_SJ^p)K^\times$, 
$$ \psi_3(\psi_2( u_{S \setminus T} j_1^p t)) =\psi_3(\psi_2( \gamma))= \gamma^{-1}\tilde{u}_S \psi_3(\psi_2( \gamma)) \in J^p$$
so $\psi_3\circ\psi_2$ is trivial and the sequence is a complex.

\smallskip

For  exactness at the second term, consider $x \in \ker(\psi_2)$. 
We immediately see
$x$ is a $p$th power at all primes of $T$ and since $x \in U_{S\setminus T} J^p  \prod_{\q \in T} K^\times_\q$ it is a $p$th power at all primes of $S\setminus T$, so $x$ is a $p$th power at all primes of $S$, and away from $S$, up to $p$th powers, it is a unit. Thus $ x \in U_SJ^p$. By hypothesis $x \in J^pK^\times$ so $x \in  \mbox{im}(\psi_1)$.

\smallskip
We now check exactness at the third term. Let $x \in \ker(\psi_3)$. 
Then $\psi_3(x) =  u_S j^p \gamma$ where $u_S \in U_S$, $j$ is an idele and $\gamma \in K^\times$.
As $\psi_3$ is `extension from $T$ to the ideles by $1$', 
we see that for  $\q \in S\setminus T$, the $\q$-component of $u_Sj^p\gamma$ is $1$.
  As $u_S$ has component $1$ at $\q \in S$ and is a unit outside of $S$,
we see $j^p\gamma$ has  $\q$-component $1$ at $\q \in  S \setminus T$ and is a unit outside of $S$.
This implies $$j^p \gamma \in (U_{S\setminus T} \prod_{\q \in T} K^\times_\q)\cap J^pK^\times  \subseteq
(U_{S\setminus T} J^p  \prod_{\q \in T} K^\times_\q) \cap J^pK^\times$$
and $\psi_2(j^p\gamma)$ has the same $T$-components as $x$ so $x \in \mbox{im}(\psi_2)$
as desired. The sequence at the beginning of the proof is exact.

\smallskip

Recall that for a group $G$, we write $G^{p,el}$ for its maximal abelian quotient  that is a vector space over $\F_p$. 
Using Lemma~\ref{lemm:Shaprops} $(vi)$ 
with $T$ replaced by $S$ gives the first term of (\ref{eq:exact}) below. 
The second term comes from Lemma~\ref{lemm:Shaprops} $(vii)$ with 
$S$ playing the same role and~$T$ there replaced by $S\setminus T$ here. 
The fourth and fifth terms come from the global Artin map.
\begin{equation}\label{eq:exact}
0 \to \frac{V_S}{K^{\times p}} \to \frac{V_{S \setminus T}^S}{K^{\times p}} \to  \prod_{\q \in T} \frac{K^\times_\q}{K^{\times p}_\q} \to G^{p,el}_S \to (G^T_S)^{p,el} \to 0.
\end{equation}
Dualizing and using  the local Artin map for the third term completes the proof.
 \end{proof}

 \subsubsection{}

We will use Proposition~\ref{proposition_10.7.9}  to control the image of 
$\Res_{N,R}:H^1(G_{N\cup R}) \to \prod_{\q \in N} H^1(G_\q)$
 when $R=\{\tilde{\mathfrak q}\}$.
 The hypothesis $\mu_p \not \subset K$ is crucial for linear disjointness of various field extensions of $K$.

 \begin{prop} \label{proposition_lifting}
  Let  $N$ be a fixed finite set of tame primes 
   with $K \subseteq L \subseteq K_N$ and let $\displaystyle{(f_\q) \in \prod_{\q \in N} H^1(G_\q,\Z/p)}$. Assume $(f_\q) $ is {\it not} in the image of the restriction map $H^1(G_{N}) \to \prod_{\q \in N} H^1(G_\q)$.
  Then there exist infinitely many finite primes $\tilde{\mathfrak q}$  such that 
  $(f_\q) \in  Im(\Res_{N,\{ \tilde{\mathfrak q}  \}})$.
  Moreover, the primes $ \tilde{\mathfrak q}$ can be chosen such that
  \begin{enumerate}
   \item[$(i)$]  $ \tilde{\mathfrak q}$ splits completely in $L/K$,
\item[$(ii)$] the  $p$-adic valuation of $N( \tilde{\mathfrak q})-1$  is larger than some given integer $m$. 
\end{enumerate}
 \end{prop}

 \begin{proof} 
 Recall first the following commutative diagram  obtained from applying Proposition~\ref{proposition_10.7.9} 
 with $(S,T)=(N,N)$ and  $(S,T)=(N \cup R, N)$. The first two vertical maps are inflation maps and thus injective. The third map is the identity. The fourth and fifth maps are duals of natural inclusions of Lemma~\ref{lemm:Shaprops} $(v)$ and $(ii)$. 
 $$\xymatrix{ 
 0 \ar@{->}[r]  & H^1(G^N_N) \ar@{->}[r]^{\Inf}  \ar@{->}[d] & H^1(G_{N}) \ar@{->}[r]^{\Res_N}  \ar@{->}[d] 
 & \prod_{\q \in N} H^1(G_\q) \ar@{->}[r]^-{\varphi}  \ar@{->}[d] &
 \CyB^N_\emptyset \ar@{->}[r]  \ar@{->}[d] & 
 \CyB_N \ar@{->}[r]  \ar@{->}[d] & 0\\
 0 \ar@{->}[r]  & H^1(G^N_{N\cup R}) \ar@{->}[r]^{\Inf} & H^1(G_{N\cup R}) \ar@{->}[r]^{\Res_{N,R}} & \prod_{\q \in N} H^1(G_\q) \ar@{->}[r]^-{\varphi_R} &
 \CyB^{N\cup R}_R \ar@{->}[r] & 
 \CyB_{N \cup R} \ar@{->}[r] & 0
 }$$
By Lemma~\ref{lemm:Shaprops}, $V^{N\cup R}_R \subset V^N_\emptyset$ and we have the field diagram below with  three of the Galois groups indicated.

$$
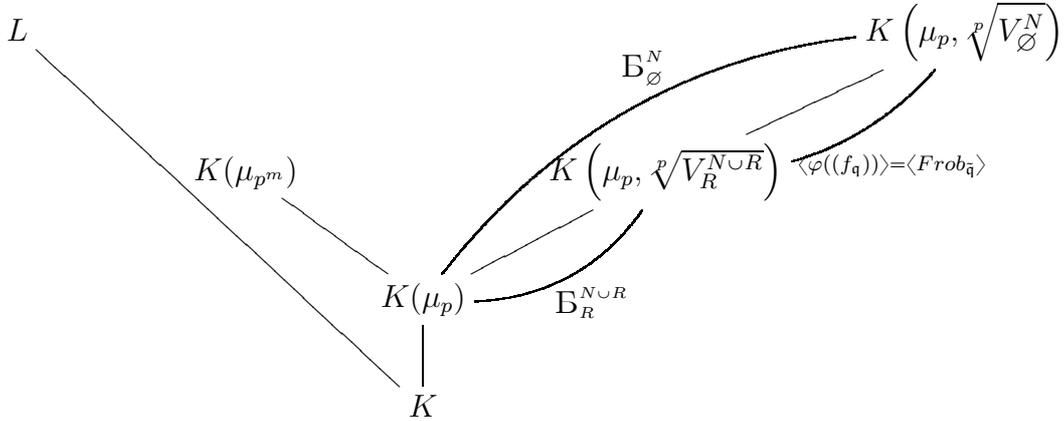
\xymatrix{ 
L& &
& & & K\left(\mu_p,\sqrt[p]{V^N_\emptyset}\right) \\
& & K(\mu_{p^m})   \ar@{-}[dr]&& K\left(\mu_p,\sqrt[p]{V^{N\cup R}_R}\right)\ar@{-}[ur] 
\ar@{-}@/^-1.50pc/[ur]_{\,\,\,\langle \varphi((f_\q)) \rangle = \langle Frob_{\tilde{\q}}\rangle}\\
& & & K(\mu_p) \ar@{-}[ur]  \ar@{-}@/^-1.50pc/[ur]_{\CyB^{N\cup R}_R} \ar@{-}@/^2.50pc/[uurr]^{\CyB^{N}_\emptyset}& \\
& & & K\ar@{-}[u]\ar@{-}[uuulll]& 
 }
 $$

 We are supposing  $(f_\q) \in  \prod_{\q \in N} H^1(G_\q) $ is not in the image of $\Res_N$. 
 Then $\varphi((f_\q)) \neq 0$ and 
 $(f_\q)$ is in the image of $\Res_{N,R} \iff \varphi_R((f_\q))=0$,  namely $(f_\q) \mapsto 0$ under the fourth vertical map. 

For  $\alpha \in V^{N}_\emptyset $ and $ \tilde{\mathfrak q} \notin N$ and of characteristic different from $p$, we have $\alpha \in U_{ \tilde{\mathfrak q}}K^{\times p}_{ \tilde{\mathfrak q}}$ so
$K\left(\mu_p,\sqrt[p]{V^{N}_\emptyset}\right) / K(\mu_p)$ is unramified at such a $ \tilde{\mathfrak q}$.
Using Chebotarev's Theorem, choose a prime $ \tilde{\mathfrak q}$ 
whose Frobenius in $\CyB^N_\emptyset = Gal\left(K\left(\mu_p, \sqrt[p]{V^{N}_\emptyset}\right) / K(\mu_p) \right)$
spans the same non-trivial line as $\varphi((f_\q)) \in \CyB^N_\emptyset $. 
 It is not hard to show  
 to show all primes of $K(\mu_p)$ above $ \tilde{\mathfrak q}$ have Frobenius elements in 
$Gal\left( K\left(\mu_p,\sqrt[p]{V^{N}_\emptyset}\right) / K(\mu_p) \right)$ that are nonzero scalar multiples of one another, so the choice is irrelevant.
Taking $R=\{ \tilde{\mathfrak q}\}$, we see
 $\alpha \in V^{N\cup R}_R$ is locally a $p$th power at $ \tilde{\mathfrak q}$ so 
$K\left(\mu_p,\sqrt[p]{V^{N \cup R}_R}\right) / K(\mu_p)$ is the maximal sub-extension of 
$K\left(\mu_p,\sqrt[p]{V^{N}_\emptyset}\right) / K(\mu_p)$ in which $ \tilde{\mathfrak q}$ splits completely so
$\varphi_R((f_\q))=0$ as desired.
 
 Since $\mu_p \not \subset K$ and $V^{N}_\emptyset \subset K$, $Gal\left( K\left(\mu_p,\sqrt[p]{V^{N}_\emptyset}\right) / K(\mu_p) \right)$
 is in the (necessarily non-trivial) $\omega$-component 
 under the action of $\Gal(K(\mu_p)/K)$), where $\omega$ is the cyclotomic character.
 As $\Gal(L(\mu_{p^m})/K(\mu_p))$ is in the trivial eigenspace, 
 $L(\mu_{p^{m}})$ and $ K\left(\mu_p,\sqrt[p]{V^{N}_\emptyset}\right) $ are linearly disjoint over $K(\mu_p)$ and 
 we can choose $ \tilde{\mathfrak q}$ to satisfy  $(i)$ and $(ii)$. 
 \end{proof}

\begin{rema}
 Proposition \ref{proposition_lifting} is reminiscent of Proposition 3.4 of \cite{RR2}, though that result  invokes Poitou-Tate duality. 
 \end{rema}

\section{Realization of  quotients of $\Gtame$}\label{section:two}

Let $G$ be a finitely generated pro-$p$ group.  
Let $d := d(G) := \dim H^1(G)$ denote the generator rank of $G$, that is, the minimum size of 
a generating set of $G$.

We will realize  groups $G$ as  quotients of $\Gtame$  by induction.
To do this we need a sequence of normal open subgroups $(H_n)_{n\geq 2}$ of $G$ such that: 
\begin{itemize}
\item[$(i)$] $H_2=[G,G]G^p$, and $\bigcap_n H_n=1$; 
 \item[$(ii)$] for every $n\geq 2$,  the quotient $H_{n}/H_{n+1}$ is isomorphic to  $\Z/p$, and is  generated by  the image of  an {\it inertial element} $y_n$ of $H_n$. Also see Definition \ref{xwhereweneedit}.
\end{itemize}


\subsection{The first step of the inductive process} \label{section_pelementary}

We will solve the main problem of this paper by induction. 
Proposition~\ref{lemm:basecase} below is the base case. Its proof uses the following result of Gras  \cite{gras}, Chapter V, Theorem 2.4:

\begin{lemm}\label{lemm:gras-munnier} Let $K$ be a number field, $T$ a finite set of primes of $K$, and $\mathfrak q \not \in T$ a prime ideal of $K$. 
There exists a $\Z/p$-extension $L/K$ ramified at exactly $\mathfrak q$  (i.e.~ it is ramified at $\q$ and at no other prime) 
and such that a finite set of primes $T$ splits completely in $L$ if and only if $\mathfrak q$ splits completely in $K\left(\mu_p,\sqrt[p]{V^T_\emptyset}\right)/K$.
\end{lemm}

This lemma allows us to prove the base case of the induction that will be the proof of Theorem~\ref{maintheorem}
which is also Theorem~\ref{TheoremA}.

\begin{prop}\label{lemm:basecase}  Let $K$ be a number field. 
There exists an elementary abelian $p$-extension $L/K$ with $\Gal(L/K)\cong (\Z/p)^d$ ramified at
$\{ \wp_1,\dots\wp_d\}$ such that $D_{\wp_i}\cong \Z/p$, where $D_{\wp_i}$ is the decomposition group of $\wp_i$ in $L/K$.  One can also choose the $\wp_i$ such that the  $p$-adic valuation of $N(\wp_i)-1$ is larger than some given integer $n_i$, for each $i$.
\end{prop}
\begin{proof}
We want an extension $L/K$ with $\Gal(L/K) \cong G\cong (\Z/p)^d$ where the primes $\wp_i$ that ramify  have no residue field extension.
This argument is a variant of   the split case of the Theorem of Scholz-Reichardt and uses the governing 
extension $K\left(\mu_p,\sqrt[p]{V^T_\emptyset}\right)/K$ of Lemma~\ref{lemm:gras-munnier}. 
For any $n_1\geq 1$, we first use Chebotarev's theorem to  choose $\wp_1$ of $K$ that splits completely in $K\left(\mu_{p^{n_1}},\sqrt[p]{V_{\emptyset}}\right)/K$.

By Lemma~\ref{lemm:gras-munnier}, we see there is a $\Z/p$-extension $L_1/K$ ramified at exactly $\wp_1$.

Now set $T=\{\wp_1\}$. Assume  $n_2 \geq 1$ is given and apply Lemma~\ref{lemm:gras-munnier}
with the {additional} requirements that $\wp_2$ splits completely in $L_1/K$ and $K(\mu_{p^{n_2}})/K$.
Then there is a $\Z/p$-extension $L_2/K$ ramified at  exactly $\wp_2$ in which $\wp_1$ splits. 
As $\wp_2$ splits in $L_1/K$, we see both $\wp_1$ and $\wp_2$ have decomposition group $\Z/p$ in $\Gal(L_1L_2/K) \cong (\Z/p)^2$.

Now set $T=\{\wp_1,\wp_2\}$ and use Lemma~\ref{lemm:gras-munnier}  to find a $\wp_3$ that also splits completely in
$L_1L_2/K$ and  $K(\mu_{p^{n_3}})/K$. 
Continuing in this fashion, set $L=L_1L_2\cdots L_d$ to obtain the result.
\end{proof}

 \subsection{The embedding problem}
 
 Let $G$ be a finitely generated 
 pro-$p$ group filtered by a sequence of normal subgroups $G^p[G,G]=H_2 \supset H_3 \supset \cdots$ of $G$ such that  for $n\geq 2$,  $H_n/H_{n+1}\cong \Z/p$. 

Consider the central extension
$$1\to H_n/H_{n+1}\to G/H_{n+1}\overset{g_n}\to G/H_n\to 1,$$
where  ${g_n}$ is the natural map and  $H_n/H_{n+1}\cong \Z/p$ so the $p$-group $G/H_n$ acts trivially on $H_n/H_{n+1}$.
Since $H_2=G^p[G,G]$,  the Frattini subgroup of $G$, we have for $n\geq 2$,
$$d(G) = d(G/H_n) = d(G/H_{n+1}),$$
where $d$ denotes the minimal number of generators.
This implies the group extension
is not split.

 Let $\Gamma$ be a pro-$p$ group, and 
 for some $n\geq 2$, let   $f_n : \Gamma \twoheadrightarrow G/H_{n}$ be a surjective homomorphism
 and consider  the embedding problem:
  $$\xymatrix{ & & & \Gamma \ar@{.>}[ld]_-{?\rho_{n+1}} \ar@{->>}[d]^{\rho_n} \\
1 \ar[r] &H_{n}/H_{n+1}  \ar[r] & G/H_{n+1} \ar@{->>}[r]_{g_n} &  G/H_n & & (\E_n)}$$
As $G/H_{n+1}\to G/H_n$ is not split, the homomorphism $\rho_{n+1}$ is surjective if it exists. 

\medskip

The embedding problem is controlled by 
$H^2(\Gamma):=H^2(\Gamma,\Z/p)$.
Let $\varepsilon_n$ be the element in 
$H^2(G/H_n)$ corresponding to the group extension:
\begin{eqnarray}\label{es1} 1\longrightarrow H_n/H_{n+1}=\Z/p \longrightarrow G/H_{n+1} \longrightarrow G/H_{n} \longrightarrow 1.\end{eqnarray}
As $d=d(G/H_n)=d(G/H_{n+1})$ we see $\varepsilon_n \neq 0$, that is  the exact sequence $(\ref{es1})$ does not split.

Let us recall:
\begin{theo} \label{theorem_embedding}  
Let  $\Inf:H^2(G/H_n) \to H^2(\Gamma)$ be the inflation map.
The embedding problem $(\E_n)$ has a  solution if and only if $\Inf(\varepsilon_n)=0$. 
Moreover, since $n\geq 2$, any solution  is always proper, that is $\rho_{n+1}$ is surjective.
The set of solutions (modulo equivalence) of $(\E_n)$ is a principal homogeneous space over $H^1(\Gamma)$.
\end{theo}

\begin{proof}
Proposition 3.5.9 and 3.5.11 of \cite{NSW}.
\end{proof}

\subsection{The strategy} \label{section_Q}

Given $n\geq 2$, let $K_n/K$ be a Galois extension in $\Ktame$ such that $\Gal(K_n/K) \cong G/H_n$. Set 
$\Gamma_n=\Gal(K_n/K)$.

 Let $S_n$ be 
 the finite set of 
 tame primes ramified in $L_n/K$. 
Recall $K_{S_n}$ is the maximal   extension of $K$ unramified outside $S_n$ and $G_{S_n}=\Gal(K_{S_n}/K)$. Observe that $L_n \subset K_{S_n}$. 
We assume $S_n$ contains $Z_0$ as in Lemma~\ref{lemm:Shaprops}, $(iv)$
 so $\Sha_{S_n}=0$.

  \medskip
  
By Theorem \ref{theorem_embedding} applied to $\Gamma:=G_{S_n}$, we note that if there is no local obstruction  at 
any $\q \in S_n$ to lifting $G/H_n$ to $G/H_{n+1}$, then  the embedding problem $(\E_n)$ has a solution in $K_{S_n}/K$.
 
 \medskip
 
 The question is then: {\it How do we create a situation for which there is no local obstruction for every quotient of $G$?}

 \smallskip

Our strategy is as follows: By Proposition~\ref{lemm:basecase} there is a map 
$G_{S_2} \twoheadrightarrow  G/P_2(G)$ ramified at $\{ \wp_1,\dots,\wp_d\}$. 
As $S_2 = \{\wp_1,\dots,\wp_d\} \cup Z_0$, we see $\Sha_{S_2}=0$ by Lemma~\ref{lemm:Shaprops} $(iv)$.
We will later show for each $\q \in S_2$ there are lifts of  $G_\q \to G_{S_2} \twoheadrightarrow  G/P_2(G)$
to $G_\q \to G$ so
  Theorem~\ref{theorem_embedding} and $(\ref{local_global})$   give a solution to ($\E_2$).
We then  use the $H^1$ group (and the principal homogenous space property, see Theorem~\ref{theorem_embedding})  to obtain  a new solution for $(\E_2)$ with no obstructions for ($\E_3$) at $\q \in S_2$. This requires introducing ramification at a new prime 
$\tilde{\q}$ in such a manner that 
$G_{\tilde{\q}} \to G/H_2$ can be lifted to $G_{\tilde{\q}} \to G$.
Thus there is no local obstruction   to ($\E_3$) at $\tilde{\q}$  as well so a solution exists and we can repeat the process. For this, we use Proposition \ref{proposition_lifting} with the fields $K_n$ here playing the role of $L$ there.

 \subsection{Lifting local  homomorphisms}
 We retain the notations of the previous sections.
 In particular, we suppose given a sub-extension $K_n/K$ of $K_{S_n}/K$, with Galois group $\Gamma_n \cong G/H_n$.

 Recall we have the exact sequence
 $$0 \to \Sha_{S_n} \to H^2(G_{S_n}) \rightarrow  \oplus_{\q \in S_n}H^2(G_\q)$$
 where $S_n \supseteq Z_0$ {of Lemma~\ref{lemm:Shaprops} $(iv)$} so $\Sha_{S_n}=0$. Thus the embedding problem $(\E_n)$ 
 has a solution exactly when it has a local solution  for every $\q \in S_n$.
 The question is then reduced to the lifting problem of {\it ramified quotients in $G$} of tame local groups.
 
 \medskip
 
Recall that all tame primes we consider satisfy $N(\q) \equiv 1$ mod $p$.   For these $\q$ the pro-$p$ completion of  $G_\q$ is 
$\Z_p \rtimes \Z_p$.
In this pro-$p$ completion, let $\tau_\q \in G_\q$ be a generator of the inertia, and $\sigma_\q$ be the Frobenius.
They satisfy the 
unique relation $[\sigma_\q, \tau_\q]  = \tau_\q^{N(\q)-1}$. See \cite[\S 10.1 ]{Koch}.

 \medskip
 For each $\q$ in our set of primes which may be ramified, we will give a {\it local plan}, that is a homomorphism 
 $\rho_\q: G_\q \to G$ lifting $\rho_{\q,n}: G_\q \to G/H_n$. 
   $$\xymatrix{ &  G_\q \ar@{.>}[ld]_{?\rho_\q} \ar@{->}[d]^{{\rho}_{\q,n}} \\
   G 
   \ar@{->>}[r]
   & G/H_n  }$$
 For  $\q \in Z_0$   we choose the local plan to be any unramified map from $G_\q \to G$ lifting the image of $\sigma_\q \in G/H_n$ to an element of $G$. The  existence follows immediately from the fact that there are no obstructions to lifting problems with $G=\Z_p$, namely the $p$-cohomological dimension is one. 
We explain some specific ramified local plans in \S~\ref{subsubsection:torsion} and~\ref{subsubsection:torsionfree} and
give a general overview in \S~\ref{subsubsection:local_plan}.

   \subsubsection{Torsion pro-$p$ groups}\label{subsubsection:torsion}
    
 This is the idea  of the proof of the   Scholz-Reichardt theorem.
Suppose that $G$ contains an element  $y$ of order  $p^m$. 
     Take a prime $\q$ with $p^m \mid N(\q)-1$.
Suppose a representation $\rho_{\q,n}:G_\q \longrightarrow G/H_n$ defined  by $\rho_{\q,n}(\sigma_\q)=\overline{1}$ and $\rho_{\q,n}(\tau_\q)=\bar{y}$ is given.
   Since  $y^{N(\q) -1}=1$, the map   $\rho_\q : G_\q \rightarrow G$   given by $\rho_\q(\sigma_\q)={1}$ and $\rho_\q(\tau_\q)=y$ is a   local plan; in particular,  $\rho_\q $ is a lift of $\rho_{\q,n}$ from $G/H_n$ to~$G$.
  This is   why we need to specify $v_p(N(\q) -1)$ in advance. 

    \subsubsection{Torsion-free pro-$p$ groups}\label{subsubsection:torsionfree}

   When the pro-$p$ group $G$ is torsion-free, the situation is relatively rigid. 
     
\begin{lemm}\label{lemma_quotient_uniform}
 Let $G_0$ be a torsion-free $p$-adic analytic group of dimension $d$, and let $\varphi: G_0 \twoheadrightarrow D$ be a continuous morphism of $G_0$ to an analytic group $D$ having the same dimension. Then $\varphi$ is an isomorphism.
\end{lemm}

\begin{proof}
As $\varphi$ is surjective,  $\ker(\varphi)$ is analytic of  dimension zero and hence  finite. As $G_0$ is torsion-free, $\ker(\varphi)=1$. 
\end{proof}
    
    \begin{prop} \label{prop_embedding}  Let $D$ be the decomposition group at $\q \nmid p$ in some torsion-free Galois group $G$.
    If $\rho_\q : G_\q \rightarrow D$ is such that $\rho_\q(\tau_\q)$ is infinite then $\rho_\q$ is an isomorphism.
    \end{prop}

    \begin{proof}
     The image of $\tau_\q$ being  infinite  forces the image of $\sigma_\q$ to be infinite as well.  
     As the decomposition group is $p$-adic analytic, the result follows from
   Lemma \ref{lemma_quotient_uniform}. \color{black}
    \end{proof}
    
    Proposition \ref{prop_embedding} shows that  if there is tame ramification in a torsion-free pro-$p$ group~$G$, then $G$  must ``contain $G_\q$''.
\color{black}
   
\subsection{Stably inertially generated  pro-$p$ groups} \label{section_locally_inertially}

Let $G$ be a finitely generated pro-$p$ group. Let $(P_n(G))$ denote the $p$-central series of $G$, meaning that $P_1(G) := G$, and $P_{n+1}(G) := P_n(G)^p [G,P_n(G)]$.

     \subsubsection{Some Definitions}
     
     \begin{defi} \label{xwhereweneedit} 
     Suppose $H$ is a pro-$p$ group, $1\neq y \in H$.
     We call $y $ inertial  if there exists $x\in H$ such that $[x,y]=y^{ap^k}$ with $a \in \Z_p^\times$ and  $k\ge 1$.  
     \end{defi} 
     \begin{rema}\label{Depth}
      1) A torsion element $y\neq 1$ is inertial. Indeed, take $x=1$, $a=1$ and $p^k$ to be the order of $y$.\\   
      2) For an element $y$ as above, we can replace $x$ by a suitable power $x^{p^t}$. This allows us to  assume  $x \in H_{n+1} \subset H_2$ for any $n$. This change shifts $k$ to $k+t$.
     \end{rema}

     \begin{defi}\label{definition_inertially_generated}  A pro-$p$ group $H$ is called inertially generated  if it can be generated by inertial elements $y_1,\cdots, y_d$.  
      A finitely generated pro-$p$ group $G$ is called stably inertially generated, if the $P_n(G)$ are inertially generated.
     \end{defi}

Stably inertially generated pro-$p$ groups $G$ are FAb. Indeed, for every $n$, the abelianization $P_n^{ab}(G)$ of $P_n(G)$ is generated by (the classes of) inertial elements, which are  torsion in  $P_n^{ab}(G)$. Now it is easy to see that the finiteness of the $P_n^{ab}(G)$ implies the finiteness of $U^{ab}$ for every open subgroup $U$ of $G$.

       \begin{rema} We have $\langle \tau_\q \rangle =I_\q \subset G_\q \subset \Gtame$. If the class number of $K$ is   prime to $p$, then the $I_\q$ generate $\Gtame$ and $[\sigma_\q,\tau_\q]=\tau^{N(\q)-1} _\q$ so $\Gtame$ is inertially generated. 
   We expect it is not stably inertially generated. 
   Inertially generated pro-$p$ groups that are not stably inertially generated exist.
      \end{rema}
    
    \begin{exem} Let  $G_0:= \langle t \rangle \ltimes H$ be the semi-direct product of $H=\langle a_1,\cdots, a_{p-1} \rangle \cong \Z_p^{p-1}$  by  $\langle t \rangle$ of order $p$ with the action:
     $$ta_1t^{-1}=a_2 ;  \,\, ta_2t^{-1}=a_3; \,\, \cdots  ; \,\, ta_{p-2}t^{-1}=a_{p-1} ; \,\,ta_{p-1}t^{-1}=a_{p-1}^{-1} a_{p-2}^{-1} \cdots a_1^{-1}. $$
     This is well-defined. Note $G_0/[G_0,G_0]=\langle \overline{t}, \overline{a_1}\rangle \cong \Z/p \times \Z/p$, hence $G_0$ is generated by~$t$ and~$a_1t$. A simple computation shows that $(a_1t)^p=1$, 
     so we have the relations   $t^p=(a_1t)^p=1$ and $G_0$ is inertially generated, but not FAb (the subgroup $H$ is open) and therefore not stably inertially generated.
    \end{exem}

  \begin{rema} 
Set $\gr_n:=\gr_n(G):=P_n(G)/P_{n+1}(G)$ and $ d_n:= \dim \gr_n(G)$. As $G$ is finitely generated, we have $d_n < \infty$.  Recall that the map $x\mapsto x^p$ sends $\gr_n(G) $ to $\gr_{n+1}(G)$.
     By Nakayama's Lemma, it is easy to see that if $\gr_n$ is generated by the images of  inertial elements $y_i$, then $P_n(G)$ is inertially generated (the conjuguate of a inertial element is still inertial). In particular if $P_n(G)$ is inertially generated, the map $x\mapsto x^p$  produces  inertial elements in $P_{n+1}(G)$. 
     As we will see, the power of uniform groups is that this map induces an isomorphism between $\gr_n$ and $\gr_{n+1}$ for every $n\geq 1$; so a uniform group $G$ is stably inertially generated if and only if $G$ is inertially generated.
     \end{rema}

     \begin{exem} \label{example_inertially_generated}  
  The pro-$p$ group $\SL_2^1$ is stably inertially generated.
  
      Let $G=\SL_2^1(\Z_p)=\ker\left(\SL_2(\Z_p)\rightarrow \SL_2(\Z/p)\right)$. Set $$x=\Matrix{1}{p}{0}{1},\, y=\Matrix{1}{0}{p}{1}, \,
     z=\Matrix{1+p}{p}{-p}{1-p} \in G.$$  The group $G$ is topologically generated by the elements $x,y,z$.
     
     Given a prime $\q$ with $p\mid N(\q)-1$, 
     let $\alpha \in \Z_p$ be the square root of $N(\q)$ that is $1$ mod~$p$. Set
     $$s= \Matrix{\alpha}{0}{0}{\alpha^{-1}}, \, 
     t=\Matrix{ \frac{ \alpha+\alpha^{-1}}2 }{ \frac{ -\alpha+\alpha^{-1}}2}{ -\frac{ \alpha+\alpha^{-1}}2}{ \frac{ \alpha+\alpha^{-1}}2} \in G.$$
    It is a routine computation to check the relations
    $$[s,x]=x^{N(\q)-1},\, [s^{-1},y]=y^{N(\q)-1}, \mbox{ and } [t,z]=z^{N(\q)-1}.$$ 
These are identical to  the relation of a  tame local group $G_\q$, namely  $[\sigma_\q ,\tau_\q] =\tau_\q^{N(\q)-1}$,  where $\sigma_\q$ is a lift of the  Frobenius and $\tau_\q$ a generator of the ramification. Thus we will be able to create local plans for $G=\SL_2^1(\Z_p)$.
  One also observes that for every $n$,  the subgroups $P_n(G)$  are topologically  generated by the elements
     $x^{p^n}, y^{p^n}, z^{p^n}$, which also are compatible with tame local relations.
     \end{exem}

If $G$ is stably inertially generated then there exists a
sequence of subgroups $H_{n}$ as in the beginning of Section \ref{section:two}.
     
\begin{lemm}     \label{lemmlocalplan}
If $G$ is  stably inertially generated then there exists  a sequence of normal open subgroups $(H_n)_{n\geq 2}$ of $G$ such that: 
\begin{itemize}
\item[$(i)$] $H_1=G$, $H_2=[G,G]G^p$, and $\bigcap_n H_n=1$; 
\item[$(ii)$] the quotient $H_1/H_2 \cong (\Z/p)^d$ can be generated by the image of inertial elements $y_1, \cdots, y_d$ of $G$ for which  there exists some $x_i\in H_{2}$ with $[x_i,y_i]=y^{a_ip^{k_i}}_i$ where $a_i \in \Z^\times_p$  and $k_i\geq 1$;
 \item[$(iii)$] for every $n\geq 2$,  the quotient $H_{n}/H_{n+1}$ is isomorphic to  $\Z/p$, and is  generated by   the image of   inertial element $y \in H_n$   
for which there exists some $x\in H_{n+1}$ with $[x,y]=y^{ap^k}$ where $a \in \Z^\times_p$  and $k\geq 1$.
\end{itemize}
\end{lemm}

\begin{proof} $(i)$ follows as $G$ is a finitely generated pro-$p$ group and $(ii)$ and $(iii)$ follow from the definition of stably inertially generated and
Remark~\ref{Depth}.
    \end{proof} 
     
     \subsubsection{The local plan}\label{subsubsection:local_plan}
     Let $y \in H_n\setminus H_{n+1}$ be a inertial element.  
     By Remark~\ref{Depth} we can take $x$ in $H_{n+1}$ \color{black}  such that $[x,y]=y^{ap^k}$ for some $a\in \Z_p^\times$, and $k\geq 1$.
      Take a prime $\q$ such that $N(\q) = 1+bp^k$  
      \color{black} with $b\in \Z_p$ (we are again specifying a lower bound for $v_p(N(\q)-1)$ in advance), and consider the reduction map $\rho_{\q,n}:G_{v} \rightarrow D_{\q,n} \subset G/H_n$, where $D_{\q,n}$ is a decomposition group
      at $\q$ in $G/H_n$,
       that sends $\sigma_\q $ to $\overline{1}$ and $\tau_\q$ to $\overline{y}$.

    Set  $\displaystyle \alpha = \frac{\log_p(1+bp^k)}{\log_p(1+ap^k)} \in \Z_p$. The homomorphism
     $\rho_\q:G_\q \to G$, sending $\sigma_\q \mapsto x^\alpha $ and $\tau_\q \mapsto y$ is easily seen to be a local plan for $G_\q$ into $G$.
    In particular,  $\rho_\q$ lifts $\rho_{\q,n}$. {\it Thus there is no obstruction to lift $\rho_{\q,n}$ from $D_{\q,n} $ to $G/H_{n+1}$}.

     
     \subsection{The result}
Let  $G$ be a stably inertially generated pro-$p$ group.
     Recall the various primes we have used:
     \begin{itemize}
     \item $Z_0$ is a set of primes chosen via Lemma~\ref{lemm:Shaprops} $(iv)$. This set guarantees that for all  $S_n \supset Z_0$, the $\Sha_{S_n}$ groups we consider are trivial. The local plan for any prime in $Z_0$ is unramified.
     \item  The set $\{\wp_1,\dots,\wp_d\}$    of Proposition \ref{lemm:basecase} is chosen to give a map $\GtameS2 \twoheadrightarrow G/H_2 \cong (\Z/p)^d$ where each $I_{\wp_i}=D_{\wp_i} \subset G/H_2$ is $\Z/p$  where $I_{\wp_i}$ is the inertia group. 
   From \S \ref{subsubsection:local_plan},  there is a local plan   $G_{\wp_i} \to G$ for each $i$. 
     \item The prime  $\tilde{\q}$   of Proposition \ref{proposition_lifting} is used once we have solved  ($\E_n$)  to 
    provide a  global cohomology class that solves all the local plans at primes in $S_n$.
     We will choose  $\tilde{\q}$ so that it has a local plan and   we can continue the inductive lifting process.
     \end{itemize}

     \smallskip
     
     In this section we prove:

     \begin{theo}[Theorem \ref{TheoremA}] \label{maintheorem}
      Let $G$ be a finitely and  stably inertially generated pro-$p$ group. Then there exists a Galois extension $L/K$ in $\Ktame/K$  such that  $\Gal(L/K)\cong G$. Moreover the extension $L/K$ can be taken such that the set of primes splitting completely is infinite. 
     \end{theo}
\begin{proof} 
We consider a sequence of normal open subgroup $H_n$ of $G$ as in Lemma \ref{lemmlocalplan}. The proof is by induction.
We will explain the complete splitting at the  end of the proof.
Recall $d$ is the generator rank of $G$.

\smallskip

$\bullet$ 
Since $G$ is  stably inertially generated, we see $G=\langle y_1,\cdots, y_d \rangle $ where the  $y_i$ are 
inertial and   $x_i \in P_2(G)$ satisfy the relation $[x_i,y_i] =y_i^{a_ip^{n_i}}$.

Proposition~\ref{lemm:basecase} gives the first step. Namely we obtain $K_2/K$ with $\Gal(K_2/K) \cong (\Z/p)^d$,
$K_2/K$ is ramified at $\{\wp_1,\dots,\wp_d\}$ and $\Z/p \cong D_{\wp_i} \subset \Gal(K_2/K)$.
Set $S_2=\{\wp_1,\cdots, \wp_d\} \cup Z_0$ so $\Sha_{S_2}=0$ by Lemma~\ref{lemm:Shaprops} $(iv)$.

 Let $\rho_2: \GtameS2 \rightarrow G/H_2$ be the homomorphism sending 
$\tau_{\wp_i}$ to $\overline{y_i}$ and $\sigma_{\wp_i}$ to $\overline{1}$ and recall we have a local plan for each $G_{\wp_i}$,  $i=1,\cdots, d$.

By Theorem~\ref{theorem_embedding}
there exists a $\Z/p$-extension of $K_3'/K_2$ in $K_{S_2}/K$, Galois over $K$, solving the lifting problem ($\E_2$), that is, we have a map $\GtameS2 \twoheadrightarrow G/H_3$.  

\smallskip
$\bullet$
The problem now is that the decomposition group $D_\q$ at $\q \in S_2$ in $\Gal(K_3'/K)$ may not be liftable to $G/H_4$, that is we may be off the local plan. 
However, the local plan to $G/H_3$ does exist and by
Theorem~\ref{theorem_embedding}
 differs from our local solution to ($\E_2$) by an element of $f_\q \in H^1(G_\q)$.
\smallskip
By Lemma~\ref{lemmlocalplan}, the quotient $H_2/H_3$ is generated by the image of an inertial element $y\in G$ for which there exists some $x \in H_3$ such that $[x,y]=y^{ap^k}$.

We  use Proposition \ref{proposition_lifting} with $N=S_2$ and $R=\{ \tilde{\q} \}$ with $\tilde{\q}$  
 splitting completely in $K_2/K$, $(f_\q) \in \im(\psi_R)$
and $v_p(N(\tilde{\q})-1)=k$.

Hence there exists  $g \in H^1(G_{S_2\cup R})$ with $g|_{G_\q}= f_\q \,\forall \, \q \in S_2$.
We may act on our solution to ($\E_2$) by $g$ to produce another solution for which all local obstructions at $\q \in S_2$ to lifting to $G/H_4$ vanish; it is on all local plans at $S_2$.
We denote by $K_3$ the fixed field of the resulting homomorphism $G_{S_2\cup R}\to G/H_3$.

As  we  allowed  ramification at $\tilde{\q}$, we need to form a local plan at $\tilde{\q}$ 
compatible with our solution to ($\E_2$).
We have  $D_{\tilde{\q},3}=I_{\tilde{\q},3}= \langle \overline{y} \rangle \cong H_2/H_3 \cong \Z/p $ where $D_{\tilde{\q},3} \subset G/H_3$ is the image of $G_{\tilde{\q}}$. Our local plan is $G_{ \tilde{\q}} \to G$ where 
 $\sigma_{ \tilde{\q}} \mapsto x^\alpha$ for suitable $\alpha$ as in \S~\ref{subsubsection:local_plan}
 and $\tau_{ \tilde{\q}} \mapsto y$.
Thus we have local plans for all $\q \in S_3:=S_2 \cup \{ \tilde{\q} \}$ and we are on all of them with our new solution to ($\E_2$).

\smallskip

We then continue the process by induction. Set $L=\bigcup_n K_n \subset \Ktame$. Then $\Gal(L/K) \cong G$.

\smallskip

$\bullet$ If we want  the set of primes splitting completely in $L/K$ to be  infinite, we proceed as follows.
Let us choose a prime ${\mathfrak r}_2$ that splits completely in $K_2/K$. Set $T_2=\{{\mathfrak r}_2\}$. 
The local plan for ${\mathfrak r}_2$ is the trivial homomorphism $G_{{\mathfrak r}_2} \to G$.
As previously 
there exists a $\Z/p$-extension of $K_3'/K_2$ in $K_{S_2}/K$, Galois over $K$, solving the lifting problem ($\E_2$).
Recall that the problem is that  the decomposition group $D_\q$ at $\q \in S_2$ in $K_3'/K$ may be not be lifted in $G/H_4$, and that ${\mathfrak r}_2$ may  not  split completely in $K_3'/K$ (we can assume that ${\mathfrak r}_2$ is unramified in $K_3'/K$).
Choose $f_{{\mathfrak r}_2} \in H^1(G_{{\mathfrak r}_2})$ such that acting on our solution to ($\E_2$)  gives trivial decomposition group at ${\mathfrak r}_2\in T_2$. 
We again use   Proposition \ref{proposition_lifting} 
with $N=S_2\cup T_2$ 
to find a global $g \in H^1(G_{N \cup \{\tilde{\q} \} })$ with
$g|_{G_{\q}} =f_\q$ for all $\q \in N$.
Acting by this class on our first solution to ($\E_2$) gives a solution on the local plan at all $\q \in N$ and split completely at ${\mathfrak r}_2\in T_2$.

Now take a prime ${\mathfrak r}_3 \notin T_2$ that splits completely in $K_3/K$. 
Put $T_3=T_2\cup\{ {\mathfrak r}_3\}$ and
continue the process. For all $n\geq 2$ we have
\begin{enumerate}   \item[$(i)$] $T_2 \subset T_3 \subset \cdots \subset T_n \subset \cdots $,
 \item[$(ii)$] $\# T_n =n-1$,
 \item[$(iii)$]  and for every $n,k $, the primes of $T_{n+k}$ split in $K_n/K$.
 \end{enumerate}

Set $T=\bigcup_n T_n$. Then $T$ is infinite, and each prime  of $T$ splits completely in $L/K$.
\end{proof}

\begin{rema} We use at most 
   $\log_p|G/H_n|$ tame primes to realize $G/H_n$ as Galois quotient of $\Gtame$.
\end{rema}
  
Corollary~\ref{Unstable torsion} follows immediately from Theorem~\ref{TheoremA} and Remark~\ref{Depth}. 
It remains to prove Corollaries~\ref{corollary_padic} and~\ref{SimplePluperfect} and Theorems~\ref{NoToralQuotient} and~\ref{SplittingTheorem}. See, respectively, Theorems~\ref{ForCorF} and~\ref{ForCorG}, Corollary~\ref{ForTheoremH} and Theorem~\ref{ForTheoremI}.


\section{Linear groups and Lie algebras} \label{section:three}

The group-theoretic results in this section apply to $p=2$, though our Galois theoretic applications still require $\mu_p \not \subset K$. 

\subsection{Definitions}\label{section_padic}

Let $G$ be a  finitely generated pro-$p$ group. Recall that $P_n(G)$ denotes  the  $p$-central descending series of $G$,
 $\gr_n(G)=P_n(G)/P_{n+1}(G)$, and $d_n= \dim \gr_n(G)$.

\medskip

The  pro-$p$ group $G$ is called {\it uniform} if 
\begin{itemize}
\item $G/G^{p}$ is abelian where $G^p$ is the normal closure of the subgroup generated by all
$p$th powers in $G$,
\item  and if for every $n\geq 1$ the map
\begin{eqnarray}\label{iso_uniform} \gr_n(G) \stackrel{x\mapsto x^p}{\longrightarrow} \gr_{n+1}(G) \end{eqnarray}
induces an isomorphism. 

\end{itemize}
In this case, $\gr_n(G) \cong (\Z/p)^d$ for an integer $d$ called the dimension of $G$. In particular, $d_n=d $ for every $n\geq 2$.

\begin{exem} For $k\geq 1$, set $\SL_m^k:=\ker\left(\SL_m(\Z_p) \rightarrow \SL_m(\Z/p^{k})\right)$ (modulo $2^{k+1}$ for $p=2$). Then $G=\SL_m^1$ is uniform  of dimension $m^2-1$ 
and  for $k\geq 1$, $P_k(G)=\SL_m^k$. See \cite[Chapter 5, Theorem 5.2]{Dixon}.
\end{exem}
\begin{rema}\label{remarque_locally_uniform}
By (\ref{iso_uniform}) it is immediate that an inertially generated uniform group is stably inertially generated.
\end{rema}

\begin{defi}
A pro-$p$ group   $G$ is called $p$-adic analytic if $G$ is a closed subgroup of $\GL_m(\Z_p)$ for some $m$.
\end{defi}

Uniform groups are the primary building blocks of $p$-adic analytic groups.

\begin{theo}
 A finitely generated pro-$p$ group $G$ is $p$-adic analytic if and only if it contains a uniform group $H$ as an open subgroup.
 \end{theo}

\begin{proof}
 See \cite[Interlude A  and \S 4, Corollary 4.1]{Dixon}.
\end{proof}

\subsection{A dictionary: Lie algebras and $p$-adic analytic groups} \label{section_dictionary}

We will consider both finitely generated $\Z_p$-Lie algebras,  {\it i.e.}   $L\cong \Z_p^d$ and Lie algebras over $p$-adic fields. 

\begin{defi}
 The $\Z_p$-Lie algebra $L$ is called powerful if  $[L\,L] \subset 2p L$.
\end{defi}

One has \cite[Theorem 9.10]{Dixon}:

\begin{theo}
 There is an equivalence of categories between powerful $\Z_p$-Lie algebras~$L$  and uniform groups~$G$.
\end{theo}

As usual, this is obtained via  a $\log$ map sending a uniform group $G$ to a unique (up to isomorphism) powerful $\Z_p$-Lie algebra $L(G)$ and an exponential map $\exp$ sending  a powerful Lie algebra $L$ to a unique (up to isomorphism) powerful group $G$. We set 
$L_G:=L(G)\otimes_{\Z_p} \Q_p$.

Given a uniform group $G$
with $x,y \in L(G)$, set $\alpha=\exp(x)$ and $\beta=\exp(y)$. 
The Lie bracket is defined in $L(G)$ as follows:  {$$[x\, y]:=  \log \left(\lim_n [\alpha^{p^n},\beta^{p^n}]^{p^{-2n}} \right)$$}
where for $g\in P_{m+1}(G)$, 
we set $g^{p^{-m}}$ to be the 
unique $g_0 \in G$ such that $g_0^{p^m}=g$.  
\medskip

When $G$ is only $p$-adic analytic as opposed to uniform, one chooses a uniform open subgroup $G_0$ of $G$, and sets $L_G:=L_{G_0}$. Of course, $L_G$ does not depend on the choice of~$G_0$.

\medskip

A Lie algebra over a field is called {\it perfect} if $\brackLL=L$.
We recall a well-known result useful in our arithmetic context, e.g. for FAb pro-$p$ groups. 

\begin{prop}\label{context}
     Let $G$ be a $p$-adic analytic group $G$ with Lie algebra~$L_G$. The following assertions are equivalent:
     \begin{enumerate}
      \item[$(i)$] the pro-$p$ group $G$ is FAb;
      \item[$(ii)$] the Lie algebra $L_G$ over $\Q_p$  is perfect;
      \item[$(iii)$] the abelianization $G^{ab}$ of $G$ is finite.
      \end{enumerate}
  \end{prop}

  \begin{proof}
This is classical. See for example \cite[Proposition 3.18]{HM}.
  \end{proof}

\begin{exem}
 Semisimple Lie algebras are perfect. In particular the groups $\SL_n^k$ are FAb.
\end{exem}

  

    \subsection{Toral and Pluperfect Lie algebras} \label{section_toral_pluperfect}
    
   Throughout this section, we will always assume that Lie algebras are finite-dimensional over a field $F$. 
    
    \begin{defi} \label{definition_toral}
     A Lie algebra $L$ is called toral if for every $x \in L$, the adjoint endomorphism $\Adx: y\mapsto \brackxy$ is semisimple.
    \end{defi}

    Abelian Lie algebras are toral as the adjoint $\Adx$ is the zero map for every $x\in L$. 
    
    \begin{prop}\label{prop:PropToral} Let $L$ be a toral Lie algebra.  
    There is no element $x\in L$ such that $\Adx$ has a nonzero eigenvalue $\lambda \in F$. In particular, if for every $x \in L$, the characteristic polynomial of $\Adx$ splits over $F$, then $L$ is abelian.
    \end{prop}

    \begin{proof}
     If $\Adx$ has a non-zero eigenvalue $\lambda$ with eigenvector $y$, then 
     $$\Ady^2(x)=-\Ady(\Adx(y))=-\Ady(\lambda y)=0.$$
      By the toral hypothesis $\Ady$ is semisimple, so $0=\Ady(x) = - \Ad_x(y) = -\lambda y$, 
      contrary to the assumption that $\lambda \neq 0$.

   Suppose now that $L$  is not abelian and choose $x$   {\it not}  in the center of~$L$.  As $\Adx$ is semisimple it implies that $\Adx$ has non-trivial eigenvalues $\lambda \in F$ (by hypothesis), which is impossible by the previous observation.
    \end{proof}

In particular, if $F$ is algebraically closed, toral is equivalent to abelian.

\begin{lemm}\label{simple}Every non-trivial toral Lie algebra $L$ has a non-trivial toral quotient~$M$ which is either simple or abelian.\end{lemm}

\begin{proof}If $V$ is a finite dimensional vector space, $W$ is a subspace of $V$, and $T\colon V\to V$ is semisimple and satisfies $T(W) \subset W$, then $T$ induces a semisimple 
linear transformation on $V/W$.  Applying this to $\Adx$ for $x\in L$, it follows that every quotient of a toral Lie algebra $L$ is again toral.  Every non-trivial Lie algebra has a non-trivial quotient which is either simple or abelian.\end{proof}

    \begin{defi}
     A Lie algebra $L$ is called pluperfect if every toral quotient $L/I$, is trivial.
    \end{defi}

A pluperfect Lie algebra is perfect: If $L$ is not perfect, then $L/[L,L]$ is a non-trivial abelian Lie algebra, and therefore a toral quotient of $L$, so $L$ is not pluperfect.
   
Observe that a simple Lie algebra is either toral or pluperfect but not both.  Over an algebraically closed field, it cannot be toral, so it must be pluperfect.

\begin{prop}
    A non-trivial Lie algebra $L$ is pluperfect if and only if $L$ is perfect and $L/\Rad(L)$ is a direct sum of pluperfect simple Lie algebras.
    \end{prop}
    
    \begin{proof} Let $L$ be pluperfect.
 If $[L\,L]$ is a proper ideal of $L$, then $L/[L\,L]$ is a non-trivial abelian and hence toral quotient of $L$, so $L$ is perfect.  Moreover, $L/\Rad(L)$ is semisimple quotient of $L$,
    so it can be written as a (possibly empty) direct sum $L_1\oplus\cdots\oplus L_m$ of simple Lie algebras.  Each $L_i$ is therefore a simple quotient of $L$, so $L_i$ is not  toral and must therefore be pluperfect.
    
 Conversely, suppose $L$ has a non-trivial toral quotient $L/I$.  By Lemma~\ref{simple}, we may assume that $L/I$ is either abelian or that it is simple.  In the first case,
    $L$ cannot be perfect.  So we assume that $L/I$ is simple.  The image of $\Rad(L)$ in $L/I$ is a solvable ideal, so it must be $(0)$, so $I\supset \Rad(L)$, and we may think of
    $L/I$ as a quotient of $L/\Rad(L) = L_1\oplus\cdots \oplus L_m$ {where each $L_i$ is simple}.  However, all maximal proper ideals of a semisimple Lie algebra are kernels of projection maps $L\to L_i$,
    so $L/I$ must be isomorphic to one of the $L_i$, which means that at least one of the summands of $L/\Rad(L)$ is not pluperfect.
    \end{proof}
    
    Over an algebraically closed field, therefore, a Lie algebra is pluperfect if and only if it is perfect.
    
      \subsection{Inertially generated Lie algebras}
     
   Throughout this section, we will always assume that Lie algebras $L$ are finite-dimensional over a field $F$  of characteristic $0$. 
     We  introduce the notion of inertially generated Lie algebras as an analog  of  inertially generated pro-$p$ groups.
     
     \begin{defi} \label{InertialLieElement}A nonzero element $y$ of a Lie algebra $L$ is called inertial 
     if it is an eigenvector
      with nonzero eigenvalue of the adjoint  $\Adx$ for some~$x$.
      A Lie algebra $L$ over a field $F$ is called inertially generated if there exists an $F$-basis  $\{y_1,\cdots, y_d\}$ with each $y_i$ inertial.
     \end{defi}

 An inertial element is strongly ad-nilpotent (see \cite[\S 16.1]{H}) and therefore ad-nilpotent (see \cite[\S 15.1]{H}).
     
\begin{prop}
      Any inertially generated Lie algebra $L$     is pluperfect.      \end{prop}

\begin{proof}
     Any quotient of an inertially generated Lie algebra $L$ is again inertially generated.
     If $I$ is an ideal of $L$ such that $L/I$ is toral, then by Proposition~\ref{prop:PropToral}, $L/I$ has no inertial elements, so 
      $L/I$ is trivial.  Thus, $L$ is pluperfect.
     \end{proof}
     
The converse is true in characteristic  $0$. 
To prove this, we begin with a lemma.
 
\begin{lemm}\label{lem_nilpotent} Let $L$ be a Lie algebra 
spanned by ad-nilpotent elements.  Then the span $I$ of the set of inertial elements in $L$
 is an ideal.\end{lemm}
\begin{proof}   If $z$ is any   ad-nilpotent  element of $L$, then $\Ad_z^{\dim L}=0$ and the  function
 $$y \mapsto\exp(\,\Ad_z)(x) =  \sum_{i=0}^{\dim L-1} \frac{\Ad_z^i(y)}{i!}$$
 is   a Lie algebra  automorphism of $L$  (see \cite[\S 2.3]{H}).
 As being inertial is a characteristic property of a Lie algebra, we see that $y$ inertial implies 
$$\exp(m\,\Ad_z)(y) = \sum_{i=0}^{\dim L-1} \frac{m^i\Ad_z^i(y)}{i!}$$     
is inertial.  By the linear independence of the sequences $1,m,m^2,\ldots,m^{\dim L-1}$ as $m$ ranges from $1$ to $\dim L$, it follows that each $\displaystyle \frac{\Ad^i_z(y)}{i!}$ lies  in~$I$.
Taking $i=1$ we see  $\Ad_z$ preserves~$I$.  Since the nilpotent elements span $L$, $I$ is preserved by $\Ad_z$ for all $z\in L$, so it is an ideal.
 \end{proof}
     
 \begin{theo} \label{theo_pluperfect} A pluperfect Lie algebra $L$      is inertially generated.
     \end{theo}

     \begin{proof}
     We consider first the case that $L$ is a simple Lie algebra.  Since it is pluperfect, it is not toral, so there exists $x\in L$ such that $\Adx$ is not semisimple.
     As $L$ is semisimple and $F$ is of characteristic $0$, $x$ admits a Jordan-Chevalley decomposition $x = x_s+x_n$, with $\Ad_{x_n}$ non-zero and nilpotent.
     
     By the Jacobson-Morozov theorem, there exists an injective homomorphism $i\colon \sl_2\to L$ sending $e := \begin{pmatrix}0&1\\ 0&0\end{pmatrix}$ to $x_n$. 
     As $e$ is inertial in $\sl_2$, it follows that $x_n = i(e)$ is inertial in~$L$.
     
    Let $G$ denote the algebraic subgroup of $\GL(L)$ which stabilizes the Lie bracket.  The Lie algebra of $G$ consists of derivations of $L$, and as $L$ is semisimple,
     it coincides with~$L$.  Let $G^\circ$ denote its identity component, which is a simple algebraic group with Lie algebra $L$.
     As $F$ is perfect and infinite, $G^\circ(F)$ is Zariski-dense in $G^\circ$  \cite{Rosenlicht}, so $L$, which is irreducible as a $G^\circ$-representation, is likewise irreducible over $G^\circ(F)$.  Since it contains
     at least one inertial element and inertial elements map to inertial elements under conjugation by elements of $G^\circ(F)$, $L$ is inertially generated.
     It follows immediately that for any semisimple Lie algebra, pluperfect implies inertially generated. 
     
     Now let $L$ be an arbitrary pluperfect Lie algebra, and let $M := L/\Rad(L)$.  
     As $M$ is a quotient of $L$, it is likewise pluperfect and therefore inertially generated.  
      By the Lie-Malcev theorem, there exists an embedding of Lie algebras $j\colon M\to L$ whose composition with the quotient morphism $L\to M$ gives the identity. By the comment after Definition~\ref{InertialLieElement},
       $j(M)$ is spanned by {strongly} ad-nilpotent elements.  On the other hand, $L$ is perfect, so by \cite[\S5, Th\'eor\`eme 1]{Bourbaki}, $\Rad(L)$ is nilpotent, so $L$ is spanned by ad-nilpotent elements.  By Lemma~\ref{lem_nilpotent}, the span $I$ of inertial elements is an ideal of $L$. 
      
      Now, $j(M)\subset I$, so $\Rad(L)$ maps onto $L/I$, which implies that $L/I$ is nilpotent.  If $L/I\neq 0$, it has a non-trivial 
      abelian quotient, contrary to the fact that $L$ is perfect.  Therefore, $I=L$, and $L$ is inertially generated. 
       \end{proof}

     \begin{coro} 
      Let  $L \subset \sl_n(F)$ be    simple over $F$. Suppose that there  exists   $x\in L$ such that $\Ad_x$ has a nonzero eigenvalue.  Then $L$ is  inertially generated. In particular this is the case if $\Ad_x$ has   two different eigenvalues $\lambda_1, \lambda_2 \in F$. 
     \end{coro}

     \begin{proof} By Proposition~\ref{prop:PropToral}, $L$ is pluperfect and Theorem  \ref{theo_pluperfect} shows it is inertially generated.
     For the second part, observe that if  $x_s$ is the semisimple part of $x$,  then $x_s\in L$ since $L$ is simple.  
     In particular, it is not difficult to see that $\lambda_1-\lambda_2$ is a nonzero eigenvalue of $ad_{x_s}$ in $F$. 
     \end{proof}

\begin{prop}
Let $L$ be a simple Lie algebra and $G^\circ$ the identity component of the algebraic group of automorphisms of $L$.  Then $L$ is pluperfect if and only $G^\circ$ is isotropic, i.e., of positive rank over $F$.
\end{prop}

\begin{proof}
By \cite[Corollaire 8.5]{BT}, $G^\circ$ is isotropic if and only if it has a non-trivial unipotent subgroup defined over $F$.  Any non-zero tangent vector of such a subgroup is ad-nilpotent in $L$.

Conversely, if $G^\circ$ is pluperfect, it is inertially generated, so $L$ contains a non-zero ad-nilpotent element $x$, which determines a homomorphism of algebraic groups $t\mapsto \exp(\Ad(t\,x))$
from the additive group to $G^\circ$.  Thus, $G^\circ$ contains a unipotent subgroup, so it must be isotropic.
\end{proof}

   \subsection{Examples over local fields} \label{subsection_local}
  
  Let $F$ be $\R$ or a $p$-adic field. A simple group $G_{/F}$ is anisotropic if and only its group of $F$-points is compact  \cite[\S 6.4]{Borel}. In the real case, this amounts to $L$ being a compact Lie algebra.  In the $p$-adic case, the only anisotropic, simply connected, absolutely simple algebraic groups are inner forms of $A_n$ \cite[\S 3.3.3]{Tits}.
At the Lie algebra level,   $L$ consists of the elements of reduced trace zero in a division algebra $D$ whose center is a finite extension of $\Q_p$.  We can see this explicitly.
    
 \begin{prop}   \label{DivisionAlgebra} Let $D$ be a division algebra over $\Q_p$ so $D$ is a Lie algebra with the bracket $\brackxy=xy-yx$. Let $D_0$ be the Lie subalgebra of elements of trace zero.
  The Lie algebra $D_0$ is simple and perfect, but not pluperfect.
    \end{prop}

    \begin{proof}
   As $D_0 \otimes \overline{\Q_p} \cong \sl_n$ which is simple, $D_0$ is simple.
     One then has to verify that for every $x\in D_0$, the adjoint $\Adx$ is semisimple. But, since $D_0$ is simple, if we write $\Adx=\Ad_{x_s}+\Ad_{x_n}$, with semisimple and nilpotent parts $x_s$  and $x_n$, 
     then $x_s$ and $x_n$ are in $D_0$. Since elements of $D_0$ have multiplicative inverses,   $x_n=0$, and then $x=x_s$ implying that $\Adx=\Ad_{x_s}$. 
    \end{proof}
    
    Recall that division algebras over $\Q_p$ are {classified} by the Brauer group of $\Q_p$, which is isomorphic to $\Q/\Z$.
     For $p>2$, \color{black} let $D=(a,p)$ be the (unique up to isomorphism) nonsplit quaternion algebra over $\Q_p$, where $a$ is not a square mod $p$. Let $D_0$ be the pure quaternions corresponding to quaternion elements of zero trace, that is a simple Lie algebra $L$ of dimension $3$ which is not pluperfect. Hence  the explicitly described Lie algebra is perfect but not pluperfect:
     $$\langle x,y,z \mid  \brackxy=pz, \brackxz=pa y , \brackyz=p^2x\rangle.$$  The  uniform group is described in the following example. 
    
    \begin{exem} \label{SL4Example}Recall $p$ is odd. \color{black} Let $a \in \Z$ such that $a$  is  not a square mod $p$.
     Set $\displaystyle{U=\left(\begin{array}{cc}0&p\\1&0 \end{array}\right)}$, and consider the two following matrices of $\SL_4(\Z_p)$: $$\displaystyle{A=\left(\begin{array}{cc}U&0\\0&-U \end{array}\right) 
     \mbox{ and }B=\left(\begin{array}{cc}0&aI_2\\I_2&0 \end{array}\right)}.$$

     Then $A^2=pI_4$, $B^2=aI_4$ and $AB=-BA=\displaystyle{\left(\begin{array}{cc}0&aU\\-U&0 \end{array}\right)}$. Hence, the $\Q_p$-algebra generated by $I_4, A,B, AB$ is isomorphic to the quaternion algebra $(a,p)$.

     Set $A_0=pA$, $B_0=pB$, $C_0=pAB$, and put $x=\exp(A_0), y=\exp(B_0), z=\exp(C_0) \in \SL_4(\Z_p)$. Then the sub-group $G$ of $\SL_4(\Z_p)$ generated by $x$, $y$ and $z$ is uniform  of dimension $3$ and has toral Lie algebra. 
    \end{exem}


  \subsection{Applications}
  
  \subsubsection{Uniform groups with pluperfect Lie algebras}

    We prove our second main result.
    
     \begin{theo}[Corollary \ref{SimplePluperfect}] \label{ForCorG}Every $p$-adic analytic group $G$ with  pluperfect Lie algebra $L_G$  has a uniform open subgroup $G_0$ which is quotient of~$\Gtame$.
     \end{theo}

     \begin{proof}  Let $G$ be a $p$-adic analytic group  of Lie algebra $L_G$. 
    Theorem \ref{theo_pluperfect} implies  that  
      $L=\sum_{i=1}^d \Q_p s_i$ for some inertial elements $s_i$
     Take a powerful $\Z_p$-Lie subalgebra  of~$L$ as follows: 
      multiply the elements $s_i$ by $p^k$ for large $k$ so that $L':=\sum_{i=1}^d \Z_p s_i$   is powerful.
       Set $G_0=\exp(L')$. The group $G_0$ is uniform and generated by the exponentials of the inertial elements $x_i:=p^k s_i$. 
     In the proof of Theorem \ref{theo_pluperfect} we saw that each $x_i$ is in a  $\SL_2$-triple so there exist $y_i,z_i \in L'$, after multiplying by  $p^k$ for large $k$, such that  
     $$[z_i\,x_i]=2p^mx_i;\,\, [z_i\,y_i]=-2p^my_i;\,\, [x_i\,y_i]=p^mz_i$$ for some $m\geq 1$. 
      These relations are the same as those satisfied by the matrices  
      $$A=\Matrix{0}{p^m}{0}{0}, B=\Matrix{0}{0}{p^m}{0}, \mbox{ and } C=\Matrix{p^m}{0}{0}{-p^m}$$ in $\SL_2$.
        Exponentiating the various $\SL_2$-triples in $L$, we see $G$ contains a uniform open subgroup $G_0$ that is inertially generated. 
        By Remark \ref{remarque_locally_uniform}, $G_0$ is stably inertially generated. Theorem \ref{maintheorem} gives the result.
     \end{proof}

     Sometimes the lattice to take in $L_G$ is natural, as is the case for the linear groups $\SL_m^1$.

     \begin{exem}\label{exemple_slm} Let  $\gl_m:=\gl_m(\Z_p)$ the $\Z_p$-Lie algebra of matrices with coefficients in $\Z_p$. 
     Denote by $E_{i,j}$ the elementary matrices of$\gl_m$, where $1$ has been replaced by $p$, and for $i=1, \cdots, m-1$, set $E_{i}=E_{i,i}-E_{i+1,i+1}+E_{i,i+1}-E_{i+1,i}$.
     
     Let $\sl_m \subset \gl_m$ be the sub-Lie-algebra of $\gl_m$ generated  by the $E_{i,j}$ and $E_i$ (all have trace zero). The $\Z_p$-Lie algebra  $\sl_m$
     is powerful. Set $y_{i,j}=\exp(E_{i,j})$,   $i\neq j$, and $y_i=\exp(E_i)$.
     The $y_{i,j}$ and $y_i$ generate $\SL_m^{1}(\Z_p)$.
     
      It is easy to see   that the $y_{i,j}$  and $y_i$ are inertial so  $\SL_m^1$ is inertially generated, and then  stably inertially generated by Remark  \ref{remarque_locally_uniform}. 
     \end{exem}
     
      One then obtains part of Corollary \ref{corollary_padic}. The full result is given in \S~\ref{subsection:lifting p-adic}. 
     \begin{theo}
     For every $m\geq 2$ and $k\geq 2$,  the pro-$p$ groups $\SL_m^k(\Z_p)$ are  quotients of $\Gtame$.
     \end{theo}

  \begin{proof} The $\SL_m^k(\Z_p)$ are uniform, and stably inertially generated by Example \ref{exemple_slm}.
   Apply Theorem \ref{maintheorem}.
  \end{proof}

     \subsubsection{Toral uniform extensions}

As $\Z_p$-extensions of number fields are only wildly  ramified, 
 there is no tame ramification  in an abelian uniform extension of a number field~$K$. 
This is also a consequence of the following Proposition.

     \begin{prop} \label{prop_toral_2}
      In a toral uniform extension of a number field~$K$, there is no tame ramification.
     \end{prop}
     
     \begin{proof}
     Suppose that there is   tame ramification at $\q$ in a non-trivial toral uniform extension. By  Proposition \ref{prop_embedding} this would imply a relation of the form $[\alpha,\beta]=\beta^{N(\q)-1}$ in the uniform group which produces the  relation  $[x\,y]=\log_p(N(\q)) y$ in $L(G)$ where $x=\log(\alpha)$ and $y=\log(\beta)$. Indeed, 
     some elementary $p$-adic analysis yields 
      $${[x\,y]= \log_p \left( \lim_n [\alpha^{p^n},\beta^{p^n}]^{p^{-2n}}\right)=
      \log_p \left( \lim_n \big(\beta^{p^n(N(\q))}\big)^{p^{-2n}} \right)
      =\log_p(N(\q))y,}$$ 
      where $\log_p$ is the usual $p$-adic logarithm. Since $\log_p(N(\q))\neq 0$, this contradicts Proposition \ref{prop:PropToral}.
     \end{proof}

    \begin{coro}[Theorem \ref{NoToralQuotient}] \label{ForTheoremH}
    If the Hilbert $p$-class field tower of $K$ is finite, then
    there is no non-trivial toral uniform quotient of $\Gtame$. 
    \end{coro}
    \begin{proof}  
    Let $G$ be a non-trivial toral uniform pro-$p$ group. If $G$ is quotient of $\Gtame$, \color{black} then since the $p$-class field tower of $K$ is finite, there is tame ramification in the corresponding extension contradicting 
 Proposition \ref{prop_toral_2}.
    \end{proof}
 
 In particular, 
 the group of Example~\ref{SL4Example} 
 is not quotient of $\Gtame$.   
   
Recall the Fontaine-Mazur conjecture
for tame extensions predicts 
this result   holds without the class field tower hypothesis.
  
  \smallskip  
    To conclude this section, we finish with an extension of the previous result.
    
\begin{theo}[Theorem \ref{SplittingTheorem}] \label{ForTheoremI}
Let $d_K$ to be the absolute discriminant of $K$. \\
\noindent
(1) Take $T$ such that $\alpha_T  + \sum_{v|\infty} \alpha_\q > \log \sqrt{|d_K|} $.  Then $G^{{\rm ta},T}_K$ has no non-trivial  uniform toral  quotient. \\
(2) Assume the GRH and  $T$ such that $\alpha^{\rm GRH}_T  + \sum_{v|\infty} \alpha^{\rm GRH}_\q > \log \sqrt{|d_K|} $.  Then $G^{{\rm ta},T}_K$ has no non-trivial  uniform toral  quotient.
\end{theo}

    \begin{proof} 
    We proceed by contradiction. Suppose $G^{{\rm ta},T}_K$ has a non-trivial uniform toral quotient $G=\Gal(L/K)$ with $L \subset\Ktame$. 
    Then by Proposition \ref{prop_toral_2}, $L/K$ is unramified. 
 We use Theorem 1 and Proposition 1 of   \cite{Ihara}: in an infinite unramified extension $L/K$ one 
 unconditionally has
  $\alpha_T  + \sum_{v|\infty} \alpha_\q \leq \log \sqrt{|d_K|} $, and assuming the GRH 
 one has
  $\alpha^{\rm GRH}_T  + \sum_{v|\infty} \alpha^{\rm GRH}_\q \leq \log \sqrt{|d_K|} $, contradicting our assumption.
\end{proof}


     \subsection{Lifting to the special linear group over complete local Noetherian rings}\label{subsection:lifting p-adic}

    The result below follows immediately from Theorem \ref{maintheorem}.
      The prime $p$ is odd. 

   \begin{theo}[Corollary \ref{corollary_padic}] \label{ForCorF} For any complete Noetherian local  ring $A$ with residue field $\F_p$, $\SL_m^k(A)$ is a quotient of $\Gtame$
and can correspond to a Galois extension $L/K$ in which infinitely many primes split completely. 
   \end{theo}

   \begin{proof}  We first prove the result for the ring $A= \Z_p\ldbrack T_1,\cdots, T_n\rdbrack$. 

   The proof is an extension of Example~\ref{example_inertially_generated}, with the 
   technical difficulty that we cannot use the exponential map.
   
   First by Proposition 13.29 of \cite{Dixon}, the sequence $\left(\SL_m^k(\Z_p\ldbrack T_1,\cdots, T_n\rdbrack)\right)_k$ corresponds to the $p$-central series of $\SL_m^k(\Z_p\ldbrack T_1,\cdots, T_n\rdbrack)$.  Set $\m=(p,T_1,\cdots, T_n)$ the maximal ideal of $\Z_p\ldbrack T_1,\cdots, T_n\rdbrack)$. \color{black} Hence, it suffices to prove that for each $k\ge 1$, 
$$\SL_m^k(\Z_p\ldbrack T_1,\cdots, T_n\rdbrack)/\SL_m^{k+1}(\Z_p\ldbrack T_1,\cdots, T_n\rdbrack)\cong \m^k/\m^{k+1}\otimes_{\Z} M_m^0(\Z)$$
is spanned by the images of inertial elements  of $\SL_m^k(\Z_p\ldbrack T_1,\cdots, T_n\rdbrack)$.

We consider first the case $m=2$.  If $a_0+\cdots+a_n = k$, the following relations hold in $\SL_2^k(\Z_p\ldbrack T_1,\cdots, T_n\rdbrack)$:
\begin{align*}
\begin{pmatrix}1-p^k&0\\ 0&(1-p^k)^{-1}\end{pmatrix}
\begin{pmatrix} 1&p^{a_0} T_1^{a_1}\cdots T_n^{a_n}\\ 0 &1\end{pmatrix}
&\begin{pmatrix}(1-p^k)^{-1}&0\\ 0&(1-p^k)\end{pmatrix}\\
&= \begin{pmatrix} 1&p^{a_0} T_1^{a_1}\cdots T_n^{a_n}\\ 0 &1\end{pmatrix}^{(p^k-1)^2},\\
\begin{pmatrix}(1-p^k)^{-1}&0\\ 0&1-p^k\end{pmatrix}
\begin{pmatrix} 1&0 \\ p^{a_0} T_1^{a_1}\cdots T_n^{a_n} &1\end{pmatrix}
&\begin{pmatrix}1-p^k&0\\ 0&(1-p^k)^{-1}\end{pmatrix}\\
&= \begin{pmatrix} 1&0\\ p^{a_0} T_1^{a_1}\cdots T_n^{a_n} &1\end{pmatrix}^{(p^k-1)^2}.\\
\end{align*}
Also, if 
$$N:= \begin{pmatrix} 1&1\\ -1&-1\end{pmatrix},\ 
D := \begin{pmatrix}\frac{(1-p^k)+(1-p^k)^{-1}}2& \frac{(1-p^k)-(1-p^k)^{-1}}2 \\ 
\frac{(1-p^k)-(1-p^k)^{-1}}2& \frac{(1-p^k)+(1-p^k)^{-1}}2
\end{pmatrix},
$$
then $N$ is nilpotent and $DND^{-1} = (p^k-1)^2 N$, so
$$D (I + p^{a_0} T_1^{a_1}\cdots T_n^{a_n} N) D^{-1} = (I + p^{a_0} T_1^{a_1}\cdots T_n^{a_n} N)^{(p^k-1)^2},$$
and $\m^k/\m^{k+1}\otimes_{\Z} M_2^0(\Z)$ is  spanned by inertial elements.

To finish the proof for general $m$, consider all embeddings of $\SL_2^k(\Z_p\ldbrack T_1,\cdots, T_n\rdbrack)$ in
$\SL_m^k(\Z_p\ldbrack T_1,\cdots, T_n\rdbrack)$ which come from choosing an ordered pair of standard basis vectors.
Together the images of all the inertial elements in $\SL_2^k(\Z_p\ldbrack T_1,\cdots, T_n\rdbrack)$ which we just constructed will span
$\m^k/\m^{k+1}\otimes_{\Z} M_m^0(\Z)$ because the span of all images of $M_2^0(\Z)$ in $M_m^0(\Z)$ obtained by choosing pairs of basis elements of $\Z^m$ generates $M_m^0(\Z)$.

\medskip

For the general case, observe first that  $A$ is isomorphic to  $\Z_p\ldbrack T_1,\cdots, T_n\rdbrack/I$ for some ideal $I$. It is then sufficient to prove that the reduction map $\SL_m^k(\Z_p\ldbrack T_1,\cdots, T_n\rdbrack) \rightarrow \SL_m^k(A)$ is surjective.
Then take $x \in \SL_m^k(A)$ and lift it  to an $m\times m$  matrix $c$ with entries in $\Z_p\ldbrack T_1,\cdots, T_n\rdbrack$ which is congruent to $1$ mod $\m^k$.  Say the determinant is~$d$.  Then $d$ is $1$ mod $\m^k\cap I$ and is an unit of $\Z_p\ldbrack T_1,\cdots, T_n\rdbrack$; in particular $d^{-1} $ is $1$ mod $\m^k\cap I$. Multiply the first row of $c$ by $d^{-1}$ to get $c'$. Then  $c' \in \SL_m^k(\Z_p\ldbrack T_1,\cdots, T_n\rdbrack)$ and  $c'$ mod $I$ is exactly~$x$.
   \end{proof}

The interest of Theorem \ref{ForCorF} is that the groups $\SL_m^k(\Z_p\ldbrack T_1,\cdots, T_n\rdbrack)$
provide examples of  quotients of $\Gtame$ which
are ``between''  $p$-adic Lie groups and Golod-Shafarevich groups.  There is no known infinite  pro-$p$ quotient of $\Gtame$ ramified at  finitely many primes  that is not virtually Golod-Shafarevich, i.e. contains on open subgroup that is Golod-Shafarevich. 
 See \cite[Chapter 13]{Dixon} for more details.


\end{document}